\RequirePackage{fix-cm}
\documentclass[smallextended]{svjour3}       
\smartqed  
\usepackage{amssymb}
\usepackage{graphicx}
\usepackage{color}
\usepackage{mathtools}
\usepackage{booktabs}
\usepackage{longtable}
\usepackage{enumitem}
\usepackage{float}

\newtheorem{dfn}{Definition}
\newtheorem{notation}{Notation}

\usepackage{comment}
\usepackage{multirow}
\usepackage{algorithm}
\usepackage{algpseudocode}
\usepackage{pifont}  
\usepackage{amsmath}
\usepackage{amsfonts}
\allowdisplaybreaks
\usepackage{pgfplots}
\usepgfplotslibrary{fillbetween} 
\usetikzlibrary{intersections}   
\usepackage{tikz}
\usepgfplotslibrary{groupplots}
\usetikzlibrary{positioning, arrows.meta}
\usetikzlibrary{calc}

\usepackage{pgfplotstable}
\pgfplotsset{compat=1.17}

\usepackage{url}

\usepackage{pgfplots}
\pgfplotsset{compat=1.18}

\algnewcommand\algorithmicinput{\textbf{Input:}}
\algnewcommand\Input{\item[\algorithmicinput]}
\algnewcommand\algorithmicoutput{\textbf{Output:}}
\algnewcommand\Output{\item[\algorithmicoutput]}

\usepackage{rotating}
\usepackage{pdflscape}

\usepackage{caption}
\usepackage{subfigure}

\usepgfplotslibrary{units}

\graphicspath{{figures/}} 

\allowdisplaybreaks

\usepackage{bm}

\newcommand{\negcount}[1]{\eta(#1)}

\newcommand{\rnegeigratio}[1]{\theta\left(#1\right)}
\newcommand{\thirdmoment}[1]{\mu_3\left(#1\right)}
\newcommand{\iqr}[1]{\mathrm{IQR}\left(#1\right)}
\newcommand{\outlierproportion}[1]{P_{\text{out}}\left(#1\right)}
\newcommand{\coefvar}[1]{\mathrm{CV}\left(#1\right)}

\newcommand{\hq}[2]{h_{#1}(#2)}

\DeclareMathAlphabet{\pazocal}{OMS}{zplm}{m}{n}
\newcommand{\unif}{\pazocal{U}}

\newcommand{\mineigscaled}[1]{\widetilde{\lambda}_{\min}\!\left(#1\right)}
\newcommand{\maxeigscaled}[1]{\widetilde{\lambda}_{\max}\!\left(#1\right)}

\newcommand{\Lambdamin}{\widetilde{\Lambda}_{\min}}

\newcommand{\skewx}[1]{\operatorname{skew}\!\left(#1\right)}
\newcommand{\logcondx}[1]{\operatorname{logcond}\!\left(#1\right)}

\usepackage{siunitx}

\AtBeginEnvironment{quote}{\vspace{\baselineskip}}
\AtEndEnvironment{quote}{\vspace{\baselineskip}}

\newcounter{mycount} \renewcommand{\themycount}{\arabic{mycount}}
\AtBeginEnvironment{enumerate}{\stepcounter{mycount}}
\setlist[enumerate]{label=\themycount.\arabic*,ref=\themycount.\arabic*} 
\newcommand{\revise}[1]{{\color{black}#1}}

\usepackage{hyperref}

\begin{document}


\title{Learning to \revise{Select Between LP and SDP Relaxations for} Nonconvex Quadratically Constrained Quadratic Programs%
\thanks{Part of this research was conducted while the first two authors
were master's students at Sabanc{\i} University and is partially based
on their respective master's theses \cite{mugethesis,buketthesis}.
This work was partially supported by the BAGEP Award of the Science Academy of the third author.}
}

\author{Müge Dedeoğlu\and M. Buket Özen     \and
        Burak Kocuk 
}


\institute{ 
           Müge Dedeoğlu \at
           Sabanc{\i} University, 34956 Orhanl{\i}-Tuzla-Istanbul\\
              \email{muge.dedeoglu@sabanciuniv.edu} 
            \and
            M. Buket Özen \at
              University of Edinburgh, Peter Guthrie Tait Road, King's Buildings, Edinburgh EH9 3FD \\
              \email{M.B.Ozen@sms.ed.ac.uk}           
           \and
           Burak Kocuk \at
           Sabanc{\i} University, 34956 Orhanl{\i}-Tuzla-Istanbul\\
              \email{burakkocuk@sabanciuniv.edu}     
}

\date{Received: date / Accepted: date}

\maketitle

\begin{abstract}
Quadratically constrained quadratic programs (QCQPs) are ubiquitous in optimization: Such problems arise in applications from operations research, power systems, signal processing, chemical engineering, and portfolio theory, among others. Despite their flexibility in modeling real-life situations and the recent effort to understand their properties, nonconvex QCQPs are hard to solve in practice. Most of the approaches in the literature are based on either Linear Programming (LP) or Semidefinite Programming (SDP) relaxations, each of which works very well for some problem subclasses but perform poorly on others. In this paper, we develop a relaxation selection procedure for nonconvex QCQPs that can adaptively decide whether an LP- or SDP-based approach is expected to be more beneficial by considering the instance structure. The proposed methodology relies on utilizing machine learning methods that involve features derived from spectral properties and sparsity patterns of data matrices, and once trained appropriately \revise{within a \textit{dimension-independent} setup}, the prediction model applies to any instance with an arbitrary number of variables and constraints. We develop classification and regression models under different feature-design setups, including  \revise{a dimension-independent representation, and evaluate them on synthetically generated instances and a combined benchmark dataset consisting of MINLPLib and QPLIB instances. We further investigate relaxation selection within a spatial B\&B framework, evaluating the trade-off between lower bound quality and computational time.}


\end{abstract}

\section{Introduction}
\label{intro}
\subsection{Problem Setting and Motivation}

In this paper, we study quadratically constrained quadratic programs (QCQPs), whose generic formulation is given below:
\begin{subequations} \label{eq:generic}
\begin{align}
z = \inf_{x} \ &  x^T A_0 x + 2b_0^T x + c_0 \label{eq:quadObj} \\
\text{s.t.} \ & x^T A_k x + 2b_k^T x + c_k \le 0 \quad &k&=1, \dots, m  \label{eq:quadCons} \\
\ & l \le x \le u .\label{eq:varBounds}
\end{align}
\end{subequations}
Here, \(A_k\) is a real symmetric \(n \times n\) matrix, \(b_k\) is a real \(n\)-dimensional vector, and \(c_k\) is a real scalar for \(k = 0, \dots, m\). Additionally,   \(l\) and \(u\) are real \(n\)-dimensional vectors with \(l \le u\), where the components of \(l\) and \(u\) are allowed to take values from the extended real numbers (i.e., \(\mathbb{R} \cup \{\pm \infty\}\)). A QCQP is called \textit{convex} if  $A_k$ matrices are positive semidefinite for all \(k = 0, \dots, m\), and \textit{nonconvex} otherwise.
QCQPs
are frequently encountered in various applications from science and engineering, including operations research~\cite{Benson2013}, electric power systems~\cite{Skolfield2022}, chemical engineering~\cite{Misener2009}, signal processing~\cite{Luo2010}, portfolio selection~\cite{Mencarelli2019} and graph theory~\cite{DeKlerk2010}.
Despite their flexibility in modeling real-world situations, nonconvex QCQPs are difficult to solve in practice due to their NP-hard nature~\cite{Anstreicher2009}.

The  classical approaches to solve nonconvex QCQPs  include the standard linear programming (LP), semidefinite programming (SDP) and second-order cone programming (SOCP) based hierarchies. Although these hierarchies may  converge to the global optimum under some assumptions~\cite{Sherali1995,Lasserre2001,Parrilo2003,Sherali2013,Ahmadi2019}, they typically require a high order of relaxations, therefore, they are not extensively used in practical approaches for large-scale instances. 

Instead, most studies focus on how to construct strong convex relaxations that accurately approximate the nonconvex QCQP. There are at least three ways this dual approach is used: i) By construction, relaxations provide valuable lower bounds (for minimization problems) that can be used to judge the quality of a feasible solution  obtained using a local solver or a heuristic approach. ii) Optimal solutions obtained from the convex relaxations can be used as an initial point for local solvers. iii) Convex relaxations are solved repeatedly in  spatial branch-and-bound (B\&B) algorithms used for the global optimization of nonconvex QCQPs. 
Therefore, the selection of the \textit{right} relaxation is a crucial task for nonconvex QCQPs.

Two convex relaxations that are heavily used and extensively studied in the literature are SDP relaxation \cite{Vandenberghe1996} and  LP relaxation\footnote{In this paper, unless otherwise stated, we will consider the LP relaxation obtained by the McCormick envelopes only.} \cite{McCormick1976}.
There is considerable interest in identifying sufficient conditions under which these relaxations can be tight; for example, \cite{zhang2000quadratic,ye2003new,beck2003strong,Luo2010,wang2022tightness} for studies focusing on the SDP relaxation and   \cite{qiu2023polyhedralpropertiesrltrelaxations}  for studies focusing on the LP relaxation obtained via     Reformulation-Linearization Technique (RLT). 
Although these theoretical results are insightful, the corresponding sufficient conditions are typically restrictive and do not directly resolve the relaxation-selection problem.

An interesting question arises when one would like to compare LP and SDP relaxations in terms of their relative strength. If a sufficient condition for tightness can be verified for either relaxation, then the winner is clear. However, such sufficient conditions are rarely available. Although studies   directly comparing LP and SDP relaxations in terms of their strength for general QCQPs are scarce, 
 it is interesting to note that there is well-known dichotomy across several important problem classes.  For example, for the  \revise{o}ptimal \revise{p}ower \revise{f}low \revise{p}roblem from power systems engineering, the SDP relaxation is considered as the stronger relaxation, and its superiority has also been supported theoretically in some settings \cite{kocuk2016strong}.  On the other hand, in the  \revise{p}ooling \revise{p}roblem from chemical engineering, the predominant choice is the LP relaxation in some computational studies \cite{Marandi2018}. 
It appears that these observations are not limited to these two problems: The SDP relaxation is particularly well-suited for problems such as the continuous relaxations of the \revise{s}table \revise{s}et and MAXCUT problems in graph theory \cite{goemans1995improved,Gaar2020},  \revise{o}ptimal \revise{t}ransmission \revise{s}witching \cite{fattahi2017promises,kocuk2017new} and \revise{u}nit \revise{c}ommitment problems \cite{fattahi2017conic,tuncer2022misocp} in electric power systems, as well as various problems in signal processing \cite{Luo2010}. On the other hand, other problems that tend to favor LP relaxation include     \revise{c}ircle \revise{p}acking \cite{Khajavirad2024,Taspinar2024} and \revise{l}ayout problems \cite{huchette2018strong}. More broadly, the trade-off between linear and nonlinear relaxations has also been studied in global optimization outside a machine learning framework; see, for example, \cite{Khajavirad2018,ZhangSahinidis2025,GonzalezRodriguez2025}. These studies emphasize that the choice of relaxation can strongly affect computational performance, but they do not address whether this choice can be predicted beforehand with feature engineering using supervised learning.

Our study is motivated by the observed dichotomy about known problems for which either the LP relaxation or the SDP relaxation is the predominant choice (with or without any justification). Unfortunately, for a QCQP instance without any prior knowledge, deciding the favorable relaxation is far from obvious. Solving both relaxations and selecting the stronger one is possible, but it can be computationally expensive. If one would like to avoid this, then an alternative approach should be developed. We propose using machine learning models to predict favorable relaxation based on instance characteristics. In particular, we propose machine learning models for this purpose by choosing features related to the spectral properties and sparsity patterns of data matrices $\{A_k\}_{k=0}^m$,  which are motivated by the literature mentioned above.  \revise{We call this line of research as \textit{relaxation selection (without spatial B\&B)}}.

\revise{Although predicting the stronger relaxation is useful, relaxation strength alone does not fully determine computational performance. In particular, a tighter lower bound may improve efficiency, but obtaining this bound may require considerably more computational time. This trade-off becomes especially important in spatial  B\&B algorithms for nonconvex QCQPs, where convex relaxations are repeatedly solved over progressively smaller variable domains. 
As these domains become smaller, the relative strength of LP and SDP-based relaxations may change, since their formulations respond differently to tighter variable bounds. 
Relaxations also influence the number of nodes that can be pruned and the computational effort required to explore the remaining regions. Therefore, root node comparisons alone do not provide a complete basis for evaluating relaxation selection in the context of global optimization. These considerations motivate extending our relaxation selection framework 
to a spatial B\&B setting, where we evaluate the effect of relaxation choice on global lower bound progress under computational budgets.
 We construct two independent B\&B trees for each QCQP instance, using the LP relaxation in one tree and the SDP$'$ relaxation in the other (here, SDP$'$ refers
to a strengthened SDP relaxation that exploits variable
bounds, see Section 2.2 for details). Since LP relaxations are generally less expensive to solve, they may allow more nodes to be processed within a given computational budget. In contrast, SDP$'$ relaxations may provide stronger lower bounds, potentially leading to more effective pruning. We investigate this trade-off by comparing the global lower bounds obtained from the two trees under similar computational budgets. We then train machine learning models to predict the preferred relaxation using structural features of the original QCQP. We call this line of research as \textit{relaxation selection within spatial B\&B}. In this perspective, \textit{relaxation selection without spatial B\&B} can also be seen as \textit{relaxation selection for the root node}.}

\subsection{Related Work and Research Gap}

Our work is related to the growing literature on machine learning for optimization, where learning mechanisms can be used for algorithmic decisions inside optimization models. Some studies have considered this method for tasks such as branching, decomposition, cut selection, heuristic selection, and reformulation or linearization decisions \cite{bengio2021}. In branch-and-bound settings, \cite{gasse2019} also shows that graph neural networks (GNN) can imitate strong branching and improve variable selection decisions for Mixed Integer Linear Programming (MILP) problems. Beyond branching, machine learning has also been used to guide higher-level design choices. For example, \cite{kruber2017} studies when a decomposition strategy should be applied, \cite{WEINER2023100061} investigates the prediction of decomposition quality for mixed-integer programs generated via constraint relaxation, and \cite{lodi2022} formulates the decision of whether to linearize a convex mixed-integer quadratic problem as a supervised learning task. These studies indicate that machine learning can be effectively used not only within a fixed optimization algorithm but also to predict which modeling or algorithmic choice is likely to be more beneficial for a given problem instance.

Several studies have extended this optimization and machine learning framework to nonlinear and nonconvex optimization models, which are related to our work. In quadratic optimization, \cite{baltean2019} proposes a learning-based approach for scoring positive semidefinite cutting planes, while \cite{kannan2025} develops a machine learning approximation of strong partitioning to accelerate the global optimization of nonconvex QCQPs.  In mixed-integer nonlinear programming, \cite{nannicini2011} uses support vector machines to predict probing failures, and \cite{zhang2025} studies learning-based deactivation of probing by graph convolutional networks. In polynomial optimization, \cite{ghaddar2022learning} develops a learning-based framework for spatial branching in RLT-based branch-and-bound, and \cite{gonzalez2022polynomial} uses machine learning to predict which conic constraints are most effective for strengthening RLT relaxations. These papers show the potential of learning-enhanced nonlinear optimization; however, they address different decision problems such as branching, probing, partitioning, and conic strengthening while having a fixed relaxation.

To the best of our knowledge, existing studies do not address the problem of predicting whether an LP- or SDP-based relaxation is more favorable for a given nonconvex QCQP instance. In this paper, we address this question by developing and evaluating machine learning-based instance representations for relaxation selection in nonconvex QCQPs.

\subsection{Main Contributions}

Our main methodological question in this study is how to represent a QCQP instance for relaxation selection. This decision problem can be viewed as a structured prediction task: the input is not an unstructured feature vector, but a mathematical optimization instance whose spectral, sparsity-related, and bound-related properties influence the relative strength of LP- and SDP-based relaxations. Our proposed methodology tries to encode QCQP structure through explicit, theoretically motivated descriptors derived from the optimization literature. Therefore, our goal is not to propose a new generic machine learning model but to develop and evaluate QCQP-specific instance representations for relaxation selection.


Our main contributions are as follows. First, we formulate
LP-vs.-SDP and LP-vs.-SDP$'$ relaxation selection for nonconvex
QCQPs as a supervised learning problem. Second, we develop three QCQP-specific instance
representation setups, namely fully dimension-dependent,
semi dimension-dependent, and dimension-independent,
which create an explicit accuracy-applicability trade-off.
In particular, the DI setup yields a fixed-size
representation that remains applicable across QCQPs
with different numbers of variables and constraints.
The proposed feature-based framework is interpretable,
since it is built on explicit spectral, sparsity-related,
and bound-related descriptors motivated by the
optimization literature on convex relaxations. \revise{In addition to the relaxation selection experiments, we investigate
LP-vs.-SDP$'$ relaxation selection within a spatial  B\&B  framework using benchmark instances
from MINLPLib and QPLIB libraries. We construct independent B\&B trees
using LP and SDP$'$ relaxations and compare the global lower
bounds obtained under computational budgets. These comparisons
are posed as classification and regression tasks.} We also compare our framework with a tripartite
GNN model to examine whether a model that learns
directly from the structure of a QCQP instance
can improve relaxation selection.


The remainder of this paper is structured as follows:
Section~\ref{sec:relaxations} formally introduces LP and SDP relaxations.
Section~\ref{sec:method} \revise{presents the three feature
engineering approaches, classification and regression models, and the spatial B\&B
framework used to evaluate relaxation selection. 
Section~\ref{subsec:data-generation} presents the synthetic instance generation procedure as well as the benchmark libraries used in the performance evaluation later.} 
Section~\ref{sec:results} \revise{presents the computational experiments
  for relaxation selection without and within a spatial B\&B framework.}
Finally, Section~\ref{sec:conclusions} concludes our paper with
final remarks and future research directions.

\section{Convex Relaxations for QCQPs}
\label{sec:relaxations}

We first establish some notation that will be used throughout this paper.

\begin{notation}
Let $X$ be a real symmetric  $n \times n$ matrix. Then, \( X \succeq 0 \) (resp. \( X \succ 0 \)) indicates that $X$    is positive semidefinite (resp. positive definite).
\end{notation}

\begin{notation}
The Frobenius inner product of two real symmetric \( n \times n \) matrices \( X \) and \( Y \) is denoted by \( X \bullet Y \) and is calculated as \( \sum_{i,j=1}^{n} X_{ij} Y_{ij} \).
\end{notation}

\begin{notation}
The standard basis vector in \( \mathbb{R}^n \) with 1 in the \( i \)-th position and 0 elsewhere is denoted by \( e_i \).
\end{notation}

Instead of  formulation~\eqref{eq:generic} introduced earlier, we will use a \textit{lifted formulation} by defining a matrix variable  \( X \),  which   represents the outer product \( xx^T \). This reformulation transforms the original nonconvex quadratic constraints into linear constraints in terms of   variables \( x \) and \( X \). The resulting problem is as follows:
\begin{subequations} \label{eq:genericLifted}
\begin{align}
z = \inf_{x, X} \ &  A_0 \bullet X + 2b_0^T x + c_0 \label{eq:linearizedObj} \\
\text{s.t.} \ &  A_k \bullet X + 2b_k^T x + c_k \le 0 \quad &k&=1, \dots, m \label{eq:linearizedCons} \\
\ & X = xx^T \label{eq:consistency} \\
\ & \eqref{eq:varBounds}. \notag
\end{align}
\end{subequations}
The constraint \( X = xx^T \) ensures that \( X \) is a rank-one matrix and formulation~\eqref{eq:genericLifted} is equivalent to formulation~\eqref{eq:generic}. 
Depending on how this  nonconvex constraint is relaxed, we   derive some convex relaxations of  QCQP \eqref{eq:genericLifted} below.

{\color{red}

}

\subsection{LP Relaxation}
\label{sec:LP-relax}

The LP relaxation of   the lifted QCQP \eqref{eq:genericLifted} replaces the nonconvex constraint \( X = xx^T \) with a set of linear inequalities known as \textit{McCormick envelopes}~\cite{McCormick1976}, which are used to convexify the bilinear terms \( X_{ij} = x_i x_j \) over the box of variables $(x_i,x_j)\in[l_i,u_i]\times[l_j,u_j]$:
\begin{subequations} \label{eq:genericLiftedLP}
\begin{align}
z_{\text{LP}} = \inf_{x, X} \ &  A_0 \bullet X + 2b_0^T x + c_0 \\
\text{s.t.} \ &
\eqref{eq:linearizedCons} \notag \\
\ & X_{ij} - l_j x_i - l_i x_j + l_i l_j \ge 0 &1 \le& i \le j \le n \label{eq:McCormick_ll} \\
\ & X_{ij} - u_j x_i - l_i x_j + l_i u_j \le 0 &1 \le& i \le j \le n \label{eq:McCormick_lu} \\
\ & X_{ij} - l_j x_i - u_i x_j + u_i l_j \le 0 &1 \le& i \le j \le n \label{eq:McCormick_ul} \\
\ & X_{ij} - u_j x_i - u_i x_j + u_i u_j \ge 0 &1 \le& i \le j \le n .\label{eq:McCormick_uu}
\end{align}
\end{subequations}

If a QCQP instance does not come with finite lower and upper bounds $l$ and $u$, then the LP relaxation cannot be constructed directly. 
In these cases, we assume that at least one strictly convex constraint exists, that is, the feasible region is contained in a full-dimensional ellipsoid. This assumption allows us to find an axis-aligned box that contains this ellipsoid, which can be used to derive variable bound vectors $l$ and $u$.
To be more precise, suppose that one of the constraints $x^T A x + 2b^T x + c \le 0$ defines this ellipsoid, where $A$ is a positive definite matrix. In this case, using KKT conditions, we can derive the following bounds for variable $x_i$ 
\[
l_i = -e_i^T A^{-1}b - \Delta_i \text{ and }
u_i = -e_i^T A^{-1}b + \revise{\Delta_i}
\]
where $\Delta_i = \sqrt{b^TA^{-1}b-c}\sqrt{e_i^TA^{-1}e_i}$. 
We note that the term $b^TA^{-1}b-c$ is nonnegative as otherwise the ellipsoid is empty.

\subsection{SDP Relaxation}
\label{sec:SDP-relax}

The SDP relaxation of the lifted reformulation \eqref{eq:genericLifted} replaces the nonconvex constraint \( X = xx^T \) with the convex constraint  \( X \succeq xx^T \):
\begin{subequations} \label{eq:genericLiftedSDP}
\begin{align}
z_{\text{SDP}} = \inf_{x, X} \ &  A_0 \bullet X + 2b_0^T x + c_0 \\
\text{s.t.} \ & 
 \eqref{eq:varBounds},
\eqref{eq:linearizedCons}
 \notag \\
\ &  \begin{bmatrix}
X & x \\
x^T & 1
\end{bmatrix} \succeq 0. \label{eq:psdCons} 
\end{align}
\end{subequations}
Note that using the Schur's Complement Lemma, the restriction  \( X \succeq xx^T \) is rewritten as the semidefinite representable constraint~\eqref{eq:psdCons}.

{Note that the SDP relaxation provided in \eqref{eq:genericLiftedSDP} does not utilize the variable bounds effectively. We also consider a strengthened SDP relaxation, denoted by {SDP$'$}, obtained by adding the bound-implied diagonal inequalities:
\begin{equation}
    X_{ii} \le (u_i + l_i) x_i - u_i l_i \quad \text{for } i = 1, \dots, n. \label{eq:SDP-chord}
\end{equation}
These constraints give the convex hull of   \(X_{ii} = x_i^2\) over the box $[l_i,u_i]$ (notice that constraint $X_{ii}\ge x_i^2$ is implied by constraint~\eqref{eq:psdCons}).

The general formulation of SDP$'$, obtained by incorporating these additional constraints, is as follows:
\begin{equation} 
z_{\text{SDP$'$}} =  \inf_{x, X} \  \left \{ A_0 \bullet X  + 2b_0^T x + c_0  : \     \eqref{eq:varBounds},
\eqref{eq:linearizedCons},  \eqref{eq:psdCons}, \eqref{eq:SDP-chord}  \right\}.
\end{equation}

}

\section{Methodology}
\label{sec:method}

This section outlines the supervised learning techniques utilized in this study. We develop classification and regression models to analyze how QCQP instances relate to LP and SDP relaxations. 
Instead of providing the QCQP sample in the form \eqref{eq:generic}  as input to machine learning  models, we propose using features that are related to the spectral properties and sparsity pattern of the data matrices. 
We explain how this is done  in Section \ref{subsec:feature_design}, and we introduce the supervised learning models utilized for relaxation selection  in Section \ref{subsec:label}. \revise{Finally, Section~\ref{subsec:bb_selection} presents the extension of our relaxation selection framework to the spatial B\&B setting.}

\subsection{Feature Design}\label{subsec:feature_design}

We  present three feature extraction procedures for our machine learning models that involve  characteristics of
the spectral properties and the sparsity patterns of data matrices. 
We develop three different setups called \textit{fully dimension-dependent} (fDD), \textit{semi dimension-dependent} (sDD) and \textit{dimension-independent} (DI). 
The fDD setup involves features which depend on the number of variables $n$ and constraints $m$, hence, it can only be used to predict the favorable relaxation for an instance family it is trained for.
The sDD setup alleviates the dependence on $n$ via aggregating some features, which makes this setup capable of producing predictions for different $n$ values but with fixed $m$ values.
Finally, the DI setup removes the dependence on $m$ as well by further aggregations and is able to provide predictions for any QCQP with arbitrary dimensions using a fixed number of features.
\revise{The proposed features are motivated by properties
of the QCQP data that affect the construction and strength
of convex relaxations. In particular, the eigenvalues of
the quadratic coefficient matrices provide information
about the convexity of the objective and constraints,
while sparsity patterns describe their quadratic
interaction structures. We also include information
about the availability of finite variable bounds,
which are used in constructing the relaxations.}

\subsubsection{Fully Dimension-Dependent Setup}\label{subsubsec:dim_dep}

In the first setup, the feature set is directly linked to the dimensionality of the decision variables $n$ and the number of constraints $m$ as detailed below. For instance, each eigenvalue of each data matrix is represented as a distinct feature.
Before going into details, we  introduce some notation and definitions.
\begin{notation}
Let \( A \) be a real symmetric \( n \times n \)   matrix. The vector of eigenvalues of \( A \) is denoted by \( \lambda(A) \), which is   sorted from the largest to the smallest eigenvalue. The rank of  \( A \) is denoted by $\rho(A)$. 
\end{notation}

\newcommand{\minx}[1]{\min\left(#1\right)}

\newcommand{\maxx}[1]{\max\left(#1\right)}

\newcommand{\avgx}[1]{\mathrm{avg}\left(#1\right)}

\newcommand{\mineig}[1]{\lambda_{\min}\left(#1\right)}

\newcommand{\maxeig}[1]{\lambda_{\max}\left(#1\right)}

\begin{dfn}
Let \( x \in \mathbb{R}^n \) be a vector.
 We define \(\negcount{x}\) as the number of negative values in \( x \)  as
\[
\negcount{x} = \sum_{i=1}^{n} \mathbf{1}_{\{x_i < 0\}},
\]
where \(\mathbf{1}_{\{x_i < 0\}}\) is an indicator function that equals 1 if \( x_i < 0 \) and 0 otherwise.
\end{dfn}
For example, \(\negcount{\lambda(A_k)}\) specifically yields the number of negative eigenvalues for the corresponding matrix \( A_k \).

\begin{dfn}
Let \( A \) be a real symmetric \( n \times n \)   matrix. We define the function \( \mathcal{P}_s(A) : \mathbb{R}^{n \times n} \to \{0, 1\}\)  as 
\[
\mathcal{P}_s(A) = 
\begin{cases} 
1 & \text{if the sparsity pattern of } A \text{ is } s \\
0 & \text{otherwise,}
\end{cases}
\]
where \( s \in \mathcal{S} = \{\text{Diagonal, Hollow, Bipartite, Tree, Chordal, Planar}\} \).
\end{dfn}

\revise{To identify the graph-based sparsity patterns,
we associate each matrix $A$ with an undirected graph
whose vertices correspond to the variables and whose
edges correspond to off-diagonal entries satisfying
$|A_{ij}|>10^{-8}$. The diagonal and hollow patterns
are identified directly from the matrix entries.
The sparsity patterns are not necessarily mutually
exclusive, so a matrix may satisfy multiple patterns.}
For brevity, we will use the initials of each sparsity pattern.
As an example, the function \( \mathcal{P}_D(A_k) \)   returns 1 if \( A_k \) is a diagonal matrix and 0 otherwise.

Table~\ref{tab:dim_dep} provides the detailed description of the features used in the fDD setup, categorized as spectral properties and sparsity patterns. In addition, the \textit{BoundsExist?} feature indicates whether the original problem has finite variable bounds or if the artificial bounds are derived as described in Section~\ref{sec:LP-relax}. 
 We note that the total number of features in this setup is a function of  both \(n\) and \(m\). 

\begin{longtable}{p{4cm}p{7cm}}
\caption{Description of the features in the fDD setup.} \label{tab:dim_dep} \\
\toprule
Notation & Description \\
\midrule
\multicolumn{2}{l}{\textbf{Spectral properties}} \\
\midrule
$\lambda_{j}(A_k)$ & $j$th largest eigenvalue of $A_k$, \( j = 1, \ldots, n,  k = 0, \ldots, m\). \\
\(\negcount{\lambda(A_k)}\) & Number of negative eigenvalues of $A_k$,  \( k = 0, \ldots, m \).\\
$\rho(A_k)$ & Rank of $A_k$,  \( k = 0, \ldots, m \)\\
\midrule
\multicolumn{2}{l}{\textbf{Sparsity pattern}} \\
\midrule
\(\mathcal{P}_s(A_k)\) & 1 if \(A_k\) has a sparsity pattern  \(s \in \mathcal{S}\), 0 otherwise, \(k = 0, \ldots, m\) \\

\midrule
BoundsExist? & 1 if the original problem has finite variable bounds, 0 if artificial bounds method is used \\
\midrule
Total Number of Features & \( (n + |\mathcal{S}| + 2)(m + 1) + 1\)\\
\bottomrule
\end{longtable}

\subsubsection{Semi Dimension-Dependent Setup}\label{subsubsec:n_indep}

In the second setup, the feature set is designed to be independent of the number of variables \( n \) but it is still allowed to be a function of   the number of constraints \( m \). This design strategy allows the features to encapsulate the essential characteristics of the constraints without being directly influenced by the size of the decision variable vector. For example, instead of providing all the eigenvalues, we define a ratio to represent the essential spectral properties of the constraints as given in the definition below.

\begin{dfn}
The ratio of the sum of the absolute values of negative
elements to the \(1\)-norm of a nonzero vector
\(x \in \mathbb{R}^n\) is denoted by \(\rnegeigratio{x}\), i.e.,
\[
\rnegeigratio{x} =
\frac{\sum_{i=1}^{n} |x_i| \mathbf{1}_{\{x_i < 0\}}}
{\|x\|_1}.
\]
\revise{For the zero vector, we define
\(\rnegeigratio{0}=0\).}
\end{dfn}

Table~\ref{tab:n_indep} provides the detailed description of the features used in the sDD setup.   We note that the total number of features in this setup is a function of   \(m\), but not \(n\).

\begin{longtable}{p{4cm}p{7cm}}
\caption{Description of the features in the sDD setup.} \label{tab:n_indep}\\
\toprule
Notation & Description \\
\midrule
\multicolumn{2}{l}{\textbf{Spectral properties}} \\
\midrule
$\negcount{\lambda(A_k)}\text{/} n$ & Number of negative eigenvalues of $A_k$ over $n$, \( k = 0, \ldots, m \). \\
$\rnegeigratio{\lambda(A_k)}$ & Negative eigenvalue ratio of $A_k$, \( k = 0, \ldots, m \) \\
$\mineig{A_k}$ & Smallest eigenvalue of $A_k$, \( k = 0, \ldots, m \). \\
$\maxeig{A_k}$ & Largest eigenvalue of $A_k$, \( k = 0, \ldots, m \). \\
$\rho(A_k)\text{/} n$ & Rank of $A_k$ over $n$,  \( k = 0, \ldots, m \).\\
\midrule
\multicolumn{2}{l}{\textbf{Sparsity pattern}} \\
\midrule
\(\mathcal{P}_s(A_k)\) & (same as before) \\
\midrule
BoundsExist? & (same as before) \\
\midrule
Total Number of Features &  \( ( |\mathcal{S}| + 5)(m + 1) + 1\)\\
\bottomrule
\end{longtable}

\subsubsection{Dimension-Independent Setup}
\label{subsubsec:nm_indep}

In the third setup, the feature set is designed to be independent
of both the number of variables \(n\) and the number of
constraints \(m\).
This method is particularly advantageous in two settings:
large-scale instances and when the training and test sets contain
instances of different sizes. The main idea of this approach is
that rather than providing individual information for each
constraint, we introduce aggregate statistics of all constraints.

Let us again introduce some notations.

\begin{notation}
Let
\[
\Xi,\ \Theta,\ \widetilde{\Lambda}_{\min},\
\widetilde{\Lambda}_{\max},\ \varrho \in \mathbb{R}^m
\]
denote the vectors of \revise{constraint-wise spectral and
rank properties}.
We collect them as
\[
\Pi := \{\Xi, \Theta, \widetilde{\Lambda}_{\min},
\widetilde{\Lambda}_{\max}, \varrho\}.
\]

\revise{The vectors in $\Pi$ are defined as }
\[
\Xi = \begin{bmatrix}
\negcount{\lambda(A_1)}\text{/} n & \cdots &
\negcount{\lambda(A_m)}\text{/} n
\end{bmatrix}^T,
\]
\[
\Theta = \begin{bmatrix}
\theta_{\varepsilon}(\lambda(A_1)) & \cdots &
\theta_{\varepsilon}(\lambda(A_m))
\end{bmatrix}^T,
\]
\revise{where, in the DI setup, we use a stabilized version of the negative eigenvalue ratio introduced in Section~\ref{subsubsec:n_indep}, denoted by $\theta_{\varepsilon}$:}
\[
\theta_{\varepsilon}(\lambda(A_k))
= \frac{\displaystyle\sum_i |\lambda_i(A_k)|\,
\mathbf{1}_{\{\lambda_i(A_k)<0\}}}
{\displaystyle\sum_i |\lambda_i(A_k)|+\varepsilon_{\theta}},
\qquad \varepsilon_{\theta}=10^{-6}.
\]

\revise{Since the DI setup aggregates these ratios across
constraint matrices, we introduce the stabilization parameter
$\varepsilon_\theta=10^{-6}$ to reduce numerical sensitivity
when the total absolute eigenvalue mass is close to zero.
This also ensures that the ratio remains well-defined for
zero matrices. }

\[
\widetilde{\Lambda}_{\min}
=
\begin{bmatrix}
\mineigscaled{A_1} & \cdots & \mineigscaled{A_m}
\end{bmatrix}^T,
\qquad
\widetilde{\Lambda}_{\max}
=
\begin{bmatrix}
\maxeigscaled{A_1} & \cdots & \maxeigscaled{A_m}
\end{bmatrix}^T,
\]

\[
\varrho =
\begin{bmatrix}
\rho(A_1)/n & \cdots & \rho(A_m)/n
\end{bmatrix}^T.
\]

\revise{For each sparsity pattern $s \in \mathcal{S}$,
we separately define the sparsity-indicator vector
$\Phi_s \in \mathbb{R}^m$ as}
\[
\Phi_s =
\begin{bmatrix}
\mathcal{P}_s(A_1) & \cdots & \mathcal{P}_s(A_m)
\end{bmatrix}^T.
\]

The spectrum-scaled extremal eigenvalues are defined as
\[
\mineigscaled{A_k}
= \frac{\mineig{A_k}}
{\displaystyle\sum_{i=1}^{n}|\lambda_i(A_k)|
+\varepsilon_{\mathrm{spec}}},
\qquad
\maxeigscaled{A_k}
= \frac{\maxeig{A_k}}
{\displaystyle\sum_{i=1}^{n}|\lambda_i(A_k)|
+\varepsilon_{\mathrm{spec}}},
\qquad
\varepsilon_{\mathrm{spec}}=10^{-6}.
\]

This constant prevents division by zero when the total
spectral mass
\[
\sum_{i=1}^{n}|\lambda_i(A_k)|
\]
is zero or very close to zero.

We use the spectrum-scaled versions instead of raw extremal
eigenvalues because the raw values depend on the overall
magnitude of \(A_k\), making them incomparable across
instances of different scales. The scaled versions are
dimensionless and lie in \([-1,1]\), enabling meaningful
aggregation in \(\Pi\).
\end{notation}

We will treat the vectors that provide specific information
about all constraints as ``samples'' and compute their
``statistics'' as detailed below.

\begin{notation}
We define the following statistics for a given vector \( x \in \mathbb{R}^n \):
\begin{itemize}
    \item \(\minx{x}\): The minimum value of the elements in the vector \( x \).
    \item \(\maxx{x}\): The maximum value of the elements in the vector \( x \).
    \item \(\avgx{x}\): The average (mean) value of the elements in the vector \( x \).
    \item $\thirdmoment{x}$: Third moment, calculated as $\frac{1}{n} \sum_{i=1}^{n} (x_i - \mu)^3$, where $\mu$ is the mean of $x$.
    \item Skewness, calculated as
    \[
    \skewx{x} = \frac{\frac{1}{n}\sum_{i=1}^{n}(x_i - \mu)^3}{\sigma^3 + \varepsilon_{\mathrm{skew}}},
    \qquad \varepsilon_{\mathrm{skew}} = 10^{-9},
    \]
    where $\mu$ and $\sigma$ are the mean and standard deviation of $x$. The constant \(\varepsilon_{\mathrm{skew}}\) prevents division by zero when the standard deviation is near zero. 
    \item $\iqr{x}$: Interquartile Range (IQR), calculated as $Q_3 - Q_1$, where $Q_3$ is the third quartile and $Q_1$ is the first quartile of $x$.
    \item $\outlierproportion{x}$: Proportion of outliers, defined as values more than $1.5$ times the IQR above $Q_3$ or below $Q_1$.
    \item $\coefvar{x}$: Stabilized Coefficient of Variation (CV). Since the mean \(\mu\) may be negative or near zero for eigenvalue distributions of indefinite matrices, the raw ratio \(\sigma/\mu\) is meaningless as a relative dispersion measure. Hence, we~use
    \[
    \coefvar{x} = \frac{\sigma}{|\mu| + \varepsilon_{\mathrm{CV}}},
    \qquad \varepsilon_{\mathrm{CV}} = 10^{-6},
    \]
    which is always non-negative and well-defined.
\end{itemize}
\end{notation}

As an example, \(\min(\Lambdamin)\) will help us determine the minimum of all the minimum eigenvalues of each constraint.

\begin{dfn}
Let \( \hq{q}{x} : \mathbb{R}^m \to \mathbb{R} \) be a function representing the statistic \( q \) of \( x \), where
$q\in\mathcal{Q}:=\{ \min, \max, \mathrm{avg},$
      $\textup{skewness}, \mathrm{IQR},$  $\textup{outliers}, \mathrm{CV} \}$.
\end{dfn}

For each property vector defined over the constraint matrices, we apply the statistics in \(\mathcal{Q}\) to obtain aggregate, dimension-independent descriptors. For example, \(\hq{q}{\Theta}\) provides various statistics about the negative eigenvalue ratios of the constraints.

\begin{dfn}
Let \( g_{\alpha}(x) \) be a function \(\mathbb{R}^m \to \mathbb{R}\) defined as
\[
g_{\alpha}(x) = \frac{1}{m} \sum_{i=1}^{m} \mathbf{1}_{\{x_i = \alpha\}},
\]
where \( x \in \mathbb{R}^m \) and \(\mathbf{1}_{\{x_i = \alpha\}}\) is the indicator function that equals 1 if \( x_i = \alpha \) and 0 otherwise.
\end{dfn}

This function counts the occurrences of \(\alpha\) in the vector \(x\) and normalizes by the number of constraints.

\revise{For the convexity and concavity features, we use tolerance-based spectral-mass conditions rather than exact equality of the entries of $\Theta$ to 0 or 1. Specifically, a constraint matrix $A_k$ is classified as convex when its negative eigenvalue mass does not exceed $10^{-8}$ times the sum of its total absolute eigenvalue mass and $\varepsilon_\theta=10^{-6}$. The analogous condition on positive eigenvalue mass is used to classify concave constraints. The proportions of matrices satisfying these conditions are denoted by $g_0(\Theta)$ and $g_1(\Theta)$, respectively.}

\revise{Using these classifications,} we also define
\(g_{\mathrm{overlap}}(\Theta)\) as the proportion of constraints that are
simultaneously classified as convex and concave,
\(g_{\mathrm{indef}}(\Theta)\) as the proportion of indefinite constraints
(neither convex nor concave), and \(g_{\mathrm{lin}}\) as the proportion of
constraints with negligible total spectral mass, which indicates a near-linear
quadratic term.
\revise{For the latter, the total absolute eigenvalue mass is compared with a tolerance of $10^{-10}$. These categories are not necessarily mutually exclusive.}

In addition to the aggregate statistics of \(\Pi\), we include two further feature families. The first family captures numerical conditioning: \revise{the condition number measures the variation in the magnitudes of the retained eigenvalues of a matrix.} We include the log-condition-number vector
\[
\mathcal{K} = \begin{bmatrix} \logcondx{A_1} & \cdots & \logcondx{A_m} \end{bmatrix}^T,
\]
where \(\logcondx{A_k}\) is the base-10 logarithm of the spectral condition number of \(A_k\), ignoring \revise{eigenvalues whose absolute values do not exceed} \(\varepsilon_{\mathrm{eig}} = 10^{-12}\) and stabilizing the denominator with \(\varepsilon_{\mathrm{num}} = 10^{-6}\).
\revise{The resulting condition number is restricted to $[1,10^{16}]$ and is set to $10^{16}$ when no eigenvalue exceeds the threshold.}
\revise{We then compute the aggregate statistics} \(\hq{q}{\mathcal{K}}\) for \(q \in \mathcal{Q}_{\mathrm{cond}} := \{\mathrm{avg}, \mathrm{IQR}, \mathrm{skewness}, \mathrm{CV}\}\).

The second family measures commutativity: matrices that commute share a common eigenbasis, indicating spectrally aligned constraints, while large commutator norms signal different eigenbasis and typically harder instances. We define the normalized commutator
\[
\gamma_{ij} = \frac{\|A_i A_j - A_j A_i\|_F}{\|A_i\|_F \, \|A_j\|_F + \varepsilon'}, \qquad \varepsilon' = 10^{-6},
\]
and we record
\[
\bar{\gamma} = \frac{1}{|\mathcal{P}|}\sum_{(i,j)\in\mathcal{P}} \gamma_{ij}, \qquad \hat{\gamma} = \max_{(i,j)\in\mathcal{P}} \gamma_{ij},
\]
over a selected set of matrix pairs \(\mathcal{P}\). The second family measures commutativity: matrices that commute share a common eigenbasis, indicating spectrally aligned constraints, while large commutator norms signal different eigenbasis and typically harder instances. We define the normalized commutator
\[
\gamma_{ij} = \frac{\|A_i A_j - A_j A_i\|_F}{\|A_i\|_F \, \|A_j\|_F + \varepsilon'}, \qquad \varepsilon' = 10^{-6},
\]
and we record
\[
\bar{\gamma} = \frac{1}{|\mathcal{P}|}\sum_{(i,j)\in\mathcal{P}} \gamma_{ij}, \qquad \hat{\gamma} = \max_{(i,j)\in\mathcal{P}} \gamma_{ij},
\]
over a selected set of matrix pairs \(\mathcal{P}\). \revise{The set $\mathcal{P}$ includes all objective--constraint pairs $(A_0,A_k)$, $k=1,\ldots,m$, and up to 40 constraint--constraint pairs, randomly sampled if necessary. If $\mathcal{P}$ is empty, both commutator features are set to zero. Small values of $\bar{\gamma}$ and $\hat{\gamma}$ suggest approximate simultaneous diagonalizability. However, zero matrices also yield zero commutators, so small values in instances with many linear constraints do not necessarily indicate meaningful spectral alignment.} 

Table~\ref{tab:nm_indep}  provides an overview of the notations and descriptions for the DI setup employed in our study. Different from the sDD setup introduced in Section \ref{subsubsec:n_indep}, the spectral properties   include statistics of vectors representing different eigenvalue characteristics together with proportions of convex and concave constraints. For the sparsity patterns, new ratios are defined to reflect the prevalence of diagonal, hollow, bipartite, tree, chordal and planar matrices among the constraints. 
 These new features enable us to aggregate the information in such a way that the total number of features is independent of $m$ and $n$.
 
\begin{longtable}{p{4cm}p{7cm}}
\caption{Description of the features in the DI setup.}  \label{tab:nm_indep}\\
\toprule
Notation & Description \\
\midrule
\textbf{Spectral properties} & \\
\midrule
$\negcount{\lambda(A_0)}\text{/} n$ & Number of negative eigenvalues of $A_0$ over $n$ \\
$\rnegeigratio{\lambda(A_0)}$ & Negative eigenvalue ratio of $A_{0}$ \\
$\mineigscaled{A_0}$ & Spectrum-scaled smallest eigenvalue of $A_0$ \\
$\maxeigscaled{A_0}$ & Spectrum-scaled largest eigenvalue of $A_0$ \\
$\rho(A_0)\text{/} n$ & Rank of $A_0$ over $n$\\
$g_{0}(\Theta)$ & Proportion of convex constraints\\
$g_{1}(\Theta)$ & Proportion of concave constraints\\
$g_{\mathrm{overlap}}(\Theta)$ & Proportion of constraints classified as both convex and concave \\
$g_{\mathrm{indef}}(\Theta)$ & Proportion of constraints classified as indefinite \\
$g_{\mathrm{lin}}$ & Proportion of constraints with negligible spectral mass \\
$\hq{q}{\mathcal{K}}$ & Statistic $q$ of $\mathcal{K}$, $q \in \mathcal{Q}_{\mathrm{cond}} := \{\mathrm{avg}, \mathrm{IQR}, \mathrm{skewness}, \mathrm{CV}\}$ \\ 
$\bar{\gamma}$ & Mean normalized commutator magnitude over selected matrix pairs \\
$\hat{\gamma}$ & Maximum normalized commutator magnitude over selected matrix pairs \\
{$\hq{q}{\Xi}$} & {Statistic \( q \) of \( \Xi \),  \( q \in \mathcal{Q} \)}\\
$\hq{q}{\Theta}$ & Statistic \( q \) of \( \Theta \),  \( q \in \mathcal{Q} \)\\
$\hq{q}{\widetilde{\Lambda}_{\min}}$ & Statistic $q$ of $\widetilde{\Lambda}_{\min}$, $q \in \mathcal{Q}$ \\
$\hq{q}{\widetilde{\Lambda}_{\max}}$ & Statistic $q$ of $\widetilde{\Lambda}_{\max}$, $q \in \mathcal{Q}$ \\
$\hq{q}{\varrho}$ & Statistic \( q \) of \( \varrho \),  \( q \in \mathcal{Q} \)\\
\midrule
\textbf{Sparsity pattern} & \\
\midrule
\( \mathcal{P}_s(A_0) \) & (same as before) \\
\( g_{s}(\Phi_s) \) & The ratio of matrices with sparsity pattern $s$ to the total number of constraints\\

\midrule
BoundsExist? & (same as before) \\
\midrule
Total Number of Features & $5 + 5 + |\Pi|\times|\mathcal{Q}| + 2|\mathcal{S}| + 1 + |\mathcal{Q}_{\mathrm{cond}}| + 2 = 64$ \\
\bottomrule
\end{longtable}

\subsection{Relaxation Selection as a Machine Learning Task} \label{subsec:label}

In this section, we present two supervised learning tasks, namely classification and regression, to predict the favorable relaxation among two alternatives.

In particular, suppose that we are comparing Relaxation 1 (e.g., LP) and Relaxation 2 (e.g., SDP or SDP$'$), where the former is less computationally demanding than the latter. Let us denote their objective values by \(z_1\) and \(z_2\), and define their \textit{relative difference} \(\delta\) as
\[
\delta = \frac{z_1-z_2}{|z_1|+|z_2|+\eta},
\]
where \(\eta=10^{-3}\) is a small stabilization parameter.
\revise{When the SDP relaxation is unbounded, the expression takes the indeterminate form $\infty/\infty$, which we set to 1 by convention, indicating that LP provides the better relaxation value.}

\revise{Throughout the paper, we use the term ``favorable relaxation''
with respect to relaxation strength. Relaxation 1 is considered favorable when
$\delta > -\epsilon$, including near-tie cases in which its bound is only
slightly weaker than that of Relaxation 2; in such cases, Relaxation 1 is
preferred because it is computationally less demanding. Otherwise,
Relaxation 2 is considered favorable.}

In the classification task, our aim is to predict the indicator function
\[
\mathbf{1}_{\{\delta>-\epsilon\}},
\]
which takes value 1 when Relaxation 1 gives a stronger or only slightly weaker value than Relaxation 2, where the tolerance is controlled by the parameter \(\epsilon\), and 0 otherwise.

In the regression task, our aim is to predict the target function $\delta$ directly. Notice that this value is between -1 and 1, and a negative (resp. positive)
value indicates that Relaxation 2 (resp. Relaxation 1)
\revise{provides the stronger relaxation bound}. One advantage in the regression task is that the target function $\delta$ can be interpreted as our confidence in the prediction. As a by-product of this approach, one might be willing to opt for Relaxation 1 even though the target value is slightly negative (again, recall that Relaxation 1 is computationally less demanding).

In addition to the supervised learning models based on our handcrafted features, we also evaluate a tripartite GNN for the relaxation-selection task. Our GNN operates on a tripartite graph (variable, quadratic term, constraint) and uses fixed-size node attributes for each node type (e.g., 5D for variables and quadratic terms, 2D for constraints) and also a small fixed global feature vector. The same message-passing and pooling functions are applied regardless of the number of variables $n$ and constraints 
$m$, and instance embeddings are obtained by permutation-invariant global pooling \cite {wu2025representing}. Therefore, the model’s input, output dimensions and parameters do not depend on $n$ or $m$, and this makes the approach dimension-independent and applicable across different QCQP sizes. In our setting, this graph-based model serves as a structure-aware alternative to the handcrafted DI feature setup. It is particularly relevant to our work because both approaches aim to remain applicable across varying QCQP dimensions, but they differ in how instance structure is represented: the DI setup relies on manually aggregated statistics, whereas the tripartite GNN learns directly from the underlying relational graph.





\subsection{Relaxation Selection within Spatial Branch-and-Bound}
\label{subsec:bb_selection}

\revise{We evaluate our relaxation selection framework within a spatial B\&B setting beyond the root node.
We focus on SDP$'$ rather than the standard SDP
relaxation because SDP$'$ incorporates bound-implied
inequalities, allowing it to use the variable
bounds obtained during branching. This makes SDP$'$
more suitable for the spatial B\&B framework considered
in our study. For each QCQP instance, we construct two
independent B\&B trees, using LP relaxations in one tree
and SDP$'$ relaxations in the other.
Both trees start from the same variable bounds and follow
the same node selection and branching rules. In particular,
we use breadth-first search and sequential variable selection.
Starting from the variable indexed by the current depth
modulo $n$, we select the next variable with an interval
width greater than $10^{-10}$ and divide its interval
at the midpoint.}

\revise{Before constructing the trees, we obtain a feasible
solution of the original QCQP to establish an upper bound.
For the benchmark instances, we use IPOPT to obtain a feasible
solution and verify its feasibility.
The same feasible upper bound, denoted by $UB$, is used in both
trees to ensure a consistent comparison. At each node $v$,
we solve the corresponding relaxation over the current variable
bounds to obtain a lower bound $LB_v$. The LP and SDP$'$
relaxations are updated using these bounds as the variable
domains become smaller. A node is pruned if its relaxation
is infeasible or if}
\[
LB_v \geq UB - \varepsilon_{\mathrm{prune}},
\qquad \varepsilon_{\mathrm{prune}}=10^{-8}.
\]
\revise{Otherwise, the node is divided into two child nodes
according to the branching rule, provided that a branchable
variable remains.}

\revise{We evaluate the progress of each tree using its global
lower bound. When a node is divided, its lower bound is
inherited by its children until their relaxations are solved.
At each iteration, the global lower bound is computed from
the lower bounds maintained for the regions represented
in the search tree, including inherited bounds of nodes
whose relaxations have not yet been solved.
We record the node-level relaxation values, solution times,
and global lower bounds throughout the search.}

\revise{For each QCQP instance, we generate up to 100 processed
nodes for both the LP and SDP$'$ trees.
Since solving an SDP$'$ relaxation is generally more expensive
than solving an LP relaxation, comparing the two trees after
the same number of nodes does not necessarily represent a
fair comparison in terms of computational effort.
Therefore, we also consider a time-based comparison by
estimating the computational cost of the two approaches.}
\revise{We train regression models to predict the computational
times of the LP and SDP$'$ trees using simple instance characteristics,
namely the numbers of variables, constraints, and nonzero
coefficients. The time models are calibrated using instances
for which both trees complete their prescribed node budgets,
and five-fold cross-fitting is employed to obtain time
predictions while reducing overfitting.
Let $\widehat{t}_{\mathrm{LP}}$ and
$\widehat{t}_{\mathrm{SDP}'}$ denote the predicted times
for the two relaxations. Taking 100 LP nodes as a reference,
we determine the corresponding SDP$'$ node budget~as}
\[
N_{\mathrm{SDP}'}
=
\min\left\{100,\,
\max\left\{1,\,
\left\lfloor
100\frac{\widehat{t}_{\mathrm{LP}}}
{\widehat{t}_{\mathrm{SDP}'}}
\right\rfloor
\right\}\right\}.
\]

\revise{The SDP$'$ node budget is restricted to $[1,100]$,
and the resulting comparison approximates a common
computational budget based on predicted solution times.
We then compare the global lower bounds obtained at the
corresponding node budgets and use their normalized
relative difference to construct the classification and
regression targets.
The input features are DI features introduced in Section~\ref{subsubsec:dim_dep}.}

\section{Datasets} \label{subsec:data-generation}

\revise{In this section, we explain how we generate a diverse   library of QCQPs   in  Section~\ref{subsec:data-generation-synthetic}. We also describe which instances are selected for our computational study from standard benchmark instances in Section~\ref{subsec:MINLPLib}. }

\subsection{Synthetic Instances} 
\label{subsec:data-generation-synthetic}

To evaluate the performance of the proposed relaxations, we generate a large and diverse set of synthetic QCQP instances. We generate two types of synthetic datasets denoted by \emph{Synthetic} and \emph{Synthetic$^\prime$}. 
The \emph{Synthetic} is designed to be broader, while \emph{Synthetic$^\prime$} is still diverse but designed to reflect the structural properties of the MINLPLib benchmark better. The reason for this choice is discussed in Section~\ref{subsec:MINLPLib}.

\subsubsection{Instance Generation} 

We first describe the generation of the \emph{Synthetic} dataset.
The data generation procedure, given in Algorithm~\ref{alg:generator_synth}, is designed to ensure \emph{diversity}, by covering dense, sparse, convex, and non-convex structures, 
and \emph{reproducibility}, by fixing all random seeds. Each instance $\iota=1,\dots,N$ is specified by data $( \{A^\iota_k\}_{k=0}^m,\{b^\iota_k\}_{k=0}^m,\{c^\iota_k\}_{k=0}^m,l^\iota,u^\iota)$ in the QCQP standard form. The generation procedure is designed by three binary decisions: (i) whether explicit variable bounds are present,
(ii) whether all quadratic matrices are drawn from a single family or from multiple families,
and (iii) whether the quadratic forms are bilinear or general quadratic.  
These three binary choices yield $2^3=8$ possible outcomes.
In particular, matrices $A^\iota_k$ are sampled with various spectral properties and sparsity patterns from the following families:
\begin{itemize}
    \item \textbf{Diagonal} with controlled or random eigenvalues
    \item \textbf{Random symmetric} with a prescribed number of negative eigenvalues
    \item \textbf{Positive definite}
    \item \textbf{Rank-one} (convex or concave) 
    \item \textbf{Zero-diagonal} (bilinear-like)
    \item \textbf{Graph-based sparse structures}: bipartite, tree, planar, or chordal
\end{itemize}
\revise{
When bilinear structure is selected, the generated quadratic matrices are
restricted to have zero diagonals, except for the designated positive definite
bounding constraint in instances without explicit variable bounds.
}
The eight possible outcomes are summarized in Table~\ref{tab:synthetic_outcomes}.

\begin{table}[htbp]
\centering
\small
\setlength{\tabcolsep}{4pt}
\caption{Eight generation outcomes used in the Synthetic dataset.}
\label{tab:synthetic_outcomes}
\begin{tabular}{|c|c|c|c|c|}
\hline
\textbf{Outcome} & \textbf{Bounds} & \textbf{Family Regime} & \textbf{Structure} & \textbf{Probability} \\
\hline
1 & No explicit bounds & Mixed families  & Bilinear      & 0.08 \\
2 & No explicit bounds & Mixed families  & Non-bilinear  & 0.05 \\
3 & No explicit bounds & Single family   & Bilinear      & 0.07 \\
4 & No explicit bounds & Single family   & Non-bilinear  & 0.05 \\
5 & Explicit bounds    & Mixed families  & Bilinear      & 0.20 \\
6 & Explicit bounds    & Mixed families  & Non-bilinear  & 0.10 \\
7 & Explicit bounds    & Single family   & Bilinear      & 0.25 \\
8 & Explicit bounds    & Single family   & Non-bilinear  & 0.20 \\
\hline
\end{tabular}
\end{table}

Observe that   $75\%$ of the instances have explicit box bounds, while  $25\%$ rely on artificially derived bounds. This choice places greater emphasis on bounded instances, while still preserving a significant subset of instances without explicit bounds. Since LP-based relaxations require variable bounds, bounded instances provide a more direct comparison setting, whereas the remaining instances allow us to also test the artificial-bound construction. Also, bilinear cases have slightly more weight than general quadratic ones, since they are both relevant for distinguishing the strength of LP and SDP based relaxations. 
Finally, both single-family and mixed-family regimes are included so that the dataset contains both more homogeneous and more heterogeneous quadratic structures.
The overall generation procedure is summarized in Algorithm~\ref{alg:generator_synth}.

\begin{algorithm}[h]
\caption{Synthetic Instance Generation Procedure}
\label{alg:generator_synth}
\begin{algorithmic}[1]
    \Input number of variables $n$, number of constraints $m$, number of instances $N$
    \Output QCQP instances $\{(\{A_k^\iota\},\{b_k^\iota\},\{c_k^\iota\},l^\iota,u^\iota): \iota=1,\dots,N\}$

    \For{$\iota=1,\dots,N$}
        \State Sample one of the 8 outcomes according to the prescribed probability vector
        \State Determine whether the instance has explicit bounds, uses a single or mixed family, and is bilinear or non-bilinear

        \If{single-family regime is selected}
            \State Select one matrix family $\mathcal{F}$
            \State Set $\mathcal{F}_k=\mathcal{F}$ for all $k=0,\dots,m$
        \Else
            \For{$k=0,\dots,m$}
                \State Select a matrix family $\mathcal{F}_k$
            \EndFor
        \EndIf

        \If{explicit bounds are absent}
            \State \revise{Select one constraint index $j\in\{1,\dots,m\}$ to serve as the bounding constraint}
        \EndIf

        \For{$k=0,\dots,m$}
            \If{\revise{explicit bounds are absent and $k=j$}}
                \State \revise{Generate $A_j^\iota \succ 0$}
                \State \revise{Sample $b_j^\iota \sim \mathrm{Unif}([-5,5]^n)$}
                \State \revise{Sample $c_j^\iota \sim -\mathrm{Unif}([0,10])$}
            \Else
                \State Generate $A_k^\iota$ from family $\mathcal{F}_k$
                \If{bilinear structure is selected}
                    \State Enforce zero diagonal in $A_k^\iota$
                \EndIf
                \State Sample $b_k^\iota \sim \mathrm{Unif}([-5,5]^n)$
                \State Sample $c_k^\iota \sim -\mathrm{Unif}([0,10])$
            \EndIf
        \EndFor

        \If{explicit bounds are present}
            \State Set $l^\iota=0,\;u^\iota=1$
        \Else
            \State \revise{Derive artificial bounds from the positive definite constraint indexed by $j$}
        \EndIf
    \EndFor

    \State \Return $\{(\{A_k^\iota\},\{b_k^\iota\},\{c_k^\iota\},l^\iota,u^\iota): \iota=1,\dots,N\}$
\end{algorithmic}
\end{algorithm}
\revise{The main differences between the \emph{Synthetic}
and \emph{Synthetic$^\prime$} datasets are the increased
proportion of single-family instances, the more frequent
generation of linear objective functions and constraints,
and the higher proportion of convex quadratic forms among
the remaining non-bilinear constraints. These modifications
are motivated by the structural characteristics observed
in the MINLPLib and QPLIB benchmark instances. We
regenerate the \emph{Synthetic$^\prime$} dataset to provide
training instances for the extended benchmark experiments,
while retaining the generation procedure described in
Algorithm~\ref{alg:generator_prime}.}

\begin{algorithm}[!t]
\caption{Modifications for the Synthetic$^\prime$ Dataset}
\label{alg:generator_prime}
\begin{algorithmic}[1]
\State Start from Algorithm~\ref{alg:generator_synth}
\State \revise{Increase the proportion of instances generated from a single matrix family}
\State Increase the frequency of linear objective functions by setting $A_0^\iota=0$ with high probability
\State Increase the frequency of linear non-bilinear constraints by setting some $A_k^\iota=0$, $k=1,\dots,m$
\State For the remaining non-bilinear constraints, sample convex forms more frequently than concave or indefinite ones
\State Keep all other steps unchanged
\end{algorithmic}
\end{algorithm}

We sample only negative values of $c_k^\iota$ to guarantee feasibility in both datasets. Once the QCQP instances are generated, we solve their LP relaxations with Gurobi 11.0.3, and their SDP and SDP$'$ relaxations with MOSEK 10.2.8 via the MOSEK Fusion API for Python, and record their optimal objective values. 
Additional implementation details are provided in Appendix~\ref{sec:instance_generation}.

\revise{

\subsubsection{Computational Characteristics of Synthetic Instances}
\label{sec:synthetic_characteristics}

We complement the learning experiments by examining global optimization
difficulty, relaxation tightness, and root-node bound differences for
Synthetic and Synthetic$'$ instances, with detailed results in
Appendix~\ref{app:synthetic_computational_characteristics}.

\paragraph{Global optimization difficulty and relaxation tightness.}
We draw random samples without replacement from the eligible Synthetic
and Synthetic$'$ populations to obtain representative subsets without
favoring particular computational or relaxation characteristics.
Eligibility requires valid LP and SDP$'$ root-node results, with finite
explicit or artificially derived variable bounds as appropriate.
We solve each sampled nonconvex QCQP globally using Gurobi under a fixed
time limit and optimality tolerance, with sample sizes, solver settings,
and certification outcomes reported in
Table~\ref{tab:synthetic_difficulty}.

Synthetic$'$ has a higher proportion of time-limit terminations and
larger median remaining Gurobi optimality gaps, suggesting greater
computational difficulty under the chosen settings, consistent with
its modified generation procedure. Detailed termination and gap
statistics appear in Tables~\ref{tab:synthetic_difficulty}
and~\ref{tab:synthetic_timelimit_gaps}.

For instances with certified global optima, relaxation tightness is
measured by
\begin{equation}
g_R=\frac{ z^*-L_R }{|z^*|+10^{-3}},
\label{eq:synthetic_global_gap}
\end{equation}
where $z^*$ denotes the certified global optimum and $L_R$ a finite
lower bound from relaxation $R$; smaller values indicate tighter bounds.
On the certified-optimal subset, SDP$'$ yields substantially smaller
median gaps than LP in both datasets (Table~\ref{tab:synthetic_gaps}).
These comparisons are conditional on certification within the time
limit, and the table reports the corresponding sample sizes.

\paragraph{Favorable relaxation statistics.}
We analyze which relaxation is more favorable for 320,500 Synthetic instances
across 11 generator families and 6,068 Synthetic$'$ instances without spatial B\&B. We compare LP with SDP
and SDP$'$ using $\delta_R$ from Section~\ref{subsec:label}.
A comparison is a near tie when $|\delta_R|\leq 0.01$; otherwise,
its sign identifies the stronger relaxation. Among all LP--SDP$'$ comparisons, near ties account for 37.51\% of Synthetic and 26.30\% of Synthetic$'$ instances, while SDP$'$ is strictly stronger in 48.78\% and \revise{64.82\%}, respectively. It is important to note that near ties are assigned to LP
because it is computationally less expensive.
Appendix~\ref{app:synthetic_computational_characteristics} reports
the complete root node distributions. 

}

\subsection{MINLPLib and QPLIB Instances}
\label{subsec:MINLPLib}

\revise{To evaluate our framework on benchmark instances, we
consider problems from MINLPLib \cite{minlplib} and QPLIB \cite{qplib}.
These libraries contain optimization problems with different
dimensions and structural characteristics, allowing us to
evaluate the proposed methodology on established benchmark
instances beyond synthetically generated problems.}

\revise{For MINLPLib, we initially consider instances with fewer than 300 variables belonging
to the QCQP and QCP classes. For QPLIB, we also select constrained
instances with fewer than 300 variables and relax the integrality
restrictions of integer and binary variables to obtain continuous
QCQPs.}

\revise{Since the LP and SDP$'$ relaxations require finite
bounds on variables appearing in quadratic terms, we
perform a preprocessing step to address missing variable
bounds in both MINLPLib and QPLIB instances. We first
attempt to derive finite bounds using the linear constraints
of each instance. When necessary, additional bound
construction procedures are applied. Instances for which
the required finite bounds cannot be established after
preprocessing are excluded. The remaining instances are
then parsed into the same internal QCQP representation
used for the synthetic experiments.}

\revise{For the root node experiments, we retain 169
MINLPLib instances, covering application families such as
pooling, packing and geometry, wastewater treatment,
water networks, and power systems. For the spatial B\&B
experiments, we consider a combined benchmark dataset
containing instances from both MINLPLib and QPLIB.
The QPLIB filtering procedure reduces the initial 453
instances to 168 after dimensional and constraint
restrictions, and subsequently to 143 after bound
preprocessing. After generating the spatial B\&B trees, we apply an additional filtering procedure to ensure that the select instances allow a consistent comparison between LP and SDP$'$ relaxations. We remove instances for which either relaxation encounters unresolved solver statuses or fails to produce a finite global lower bound during the search. These criteria are applied at the instance level, so both relaxation trees are excluded whenever an instance does not satisfy the requirements. This procedure results in 147 MINLPLib and 133 QPLIB instances, giving a combined benchmark dataset of 280 instances.}

\revise{For the spatial B\&B experiments,
we compare LP and SDP$'$ relaxations according to the
methodology introduced in Section~\ref{subsec:bb_selection}.
We also use the structural characteristics observed in
MINLPLib and QPLIB to guide the generation of the
Synthetic$^\prime$ dataset, as described before.}

\section{Results}
\label{sec:results}

\revise{This section presents the details of our computational experiments. We \revise{first} describe the experimental setting in Section~\ref{subsec:comp_setting}. Next, we present the results of our relaxation selection experiments on synthetically generated QCQP instances in Section~\ref{subsec:comp_results_LPvsSDP}, comparing LP with SDP and SDP$'$ relaxations. In Section~\ref{subsec:benchmark_instances}, we evaluate our framework on benchmark instances from MINLPLib and QPLIB, considering both root node relaxation selection and spatial   B\&B  experiments, including the gap--time trade-off analysis. We then investigate the contributions of individual features through feature importance analysis in Section~\ref{subsec:feature_importance}, followed by a feature ablation study in Section~\ref{subsec:feature_ablation}.}

\revise{Given the large number of experimental configurations, instance families, and feature representations considered in this study, reporting all classification and regression results in the manuscript would be impractical. Therefore, we focus on the main findings in this section and provide the complete experimental results, including detailed evaluation metrics, baseline comparisons, confusion matrices, and additional analyses, in our publicly available GitHub repository (https://github.com/mugededeoglu/learning-relaxation). To facilitate the interpretation of the reported classification accuracies, we also present the class-label distributions for the benchmark experiments in this section.}

\subsection{Computational Setting} \label{subsec:comp_setting}

\revise{To ensure reproducibility, we provide the detailed experimental configurations in Appendix~\ref{app:hyperparameters}.   These include the hyperparameters used for the classical machine learning models and the GNN, dataset splitting procedures, Synthetic$'$ generation settings, and solver configurations. The same model hyperparameters are used in relaxation selection experiment both without and within  spatial B\&B   to maintain consistency across the two settings. Following data preprocessing and feature extraction, the resulting datasets are prepared for training and evaluating the machine learning models.}
We employ supervised learning techniques for both classification and regression tasks. We include tree-based ensemble models such as Gradient Boosting, XGBoost, and Random Forest because our handcrafted QCQP descriptors produce heterogeneous data with nonlinear feature interactions. These models are known to be strong and robust empirical baselines. We also employ GNN to compare against our feature-based models.
For the classification task, we use 
GNN Classifier (GNN-C), XGB Classifier (XGB-C),
Gradient Boosting Classifier (GBC) and Random Forest Classifier (RFC).
For the regression task, we utilize GNN Regressor (GNN-R), XGB Regressor (XGB-R), Gradient Boosting Regressor (GBR) and Random Forest Regressor (RFR).
We tune those models to find the best performance using a grid search method where we specify a range of values for key hyperparameters. 
The learning experiments are implemented using Python 3.11.5 and Scikit-learn version 1.3.0. The experiments are executed on a system equipped with an Intel(R) Xeon(R) W-2145 CPU @ 3.70 GHz and 64 GB of RAM.

\subsection{Synthetic Instances}
\label{subsec:comp_results_LPvsSDP}

\revise{We first present the results of our relaxation
selection experiments on synthetically generated QCQP instances.}

\subsubsection{LP vs. SDP Comparison}
\revise{Tables~\ref{tab:classification-LPvsSDP} and
\ref{tab:regression-LPvsSDP} report the classification
and regression results, respectively, for the comparison
between LP and SDP relaxations.}
We design 17 experiments (labeled as ``ID'' in the first column of the tables) with various training and test set combinations (given under the second and third columns, respectively). In particular, 11 experiments are conducted on the Synthetic data set and 6 experiments are conducted on the Synthetic$'$ data set. Each data set is identified with a triplet $(n,m,N)$ that gives the number of variables, constraints and data points, respectively. 

For the Synthetic$'$ data set, we select six experiments to represent the main types of scenarios of interest, and we keep these same dimension combinations across the four synthetic-data result tables. Experiment 12 serves as a baseline case where all three setups are applicable. Experiment 13 is included as a more challenging setting with limited data. Experiment 14 again allows all three setups, but with a larger data set. Experiments 15 and 17 are chosen as complementary DI-only experiments to test generalization across different instance dimensions. Finally, Experiment 16 represents a more structurally complex setting with mixed training dimensions and transfer to an unseen larger dimension. Therefore, even though the number of Synthetic$'$ experiments is smaller, they still capture a representative range of scenarios.

Each experiment is run under fDD, sDD and DI setups as long as such a setup is applicable. For example, when the training and test sets have the same $n$ and $m$ (such as ID 1-5 and 12-14), all three setups are applicable. However, when the training and test sets have different $n$ but the same $m$ value (such as ID 7, 8 and 16), the fDD setup is no longer applicable. Finally, when these sets have different $m$ values (such as ID 6, 9-11, 15 and 17), only the DI setup is applicable.
In total, we have 36 feasible experiment-setup pairs, including 23 pairs for the Synthetic experiments and 13 pairs for the Synthetic$'$ experiments.
We note these experiments are designed so that the accuracy vs. applicability trade-off among different setups can be seen clearly.
%




\paragraph{Classification Model Performance}

We use the standard performance measure of \textit{accuracy} in the classification model, which is computed as follows:
Let  a vector of true labels \( y \) and a vector of predicted labels \( \hat{y} \) be given. Then,  
\begin{equation*}
\text{accuracy}(y, \hat{y}) = \frac{1}{N} \sum_{\iota=1}^{N} \mathbf{1}\left\{ y_\iota = \hat{y}_\iota  \right\},
\end{equation*}
where \( N \) is the size of the test set and \( \mathbf{1}\{\cdot\} \) is the indicator function.

We report  {accuracy}  as the main performance criterion in Classification experiments reported in Table~\ref{tab:classification-LPvsSDP} for three classifiers: Random Forest, XGB and Gradient Boosting. Our preliminary experiments have shown that Support Vector Classifier has been consistently outperformed by these three classifiers.
We also report the detailed results for the GNN-C.

For each feasible experiment-setup pair, we highlight in bold the best feature-based accuracy among XGB, GBC and RFC. Across the 36 feasible experiment setup pairs, GBC attains the largest number of best values overall, with 14 wins and participation in 2 ties. 
Nevertheless, the average accuracies of these three feature-based classifiers are very close: 0.9500, 0.9508 and 0.9474 for XGB, GBC and RFC, respectively. Therefore, while GBC remains stronger, all three classifiers perform at a consistently high level. In particular, in 32 of the 36 feasible experiment setup pairs, at least one feature-based classifier achieves an accuracy of at least 0.90. The only remaining cases are Experiment 3 under all three setups and Experiment 11 under the DI setup, which correspond to the most data-limited or structurally challenging settings. 

Let us also discuss the accuracy vs. applicability trade-off. The results do not demonstrate a uniform monotone decrease when moving from fDD to DI, which is contrary to one might expect. In the experiments where all three setups are feasible (IDs 1–5 and 12–14), the best feature-based accuracies across setups are often extremely close. Moreover, the DI setup attains the best or a tied-best feature-based accuracy in Experiments 2, 3, 5, and 13. Overall, these observations indicate that the DI setup preserves the predictive power of other setups and presents a broader applicability.

\begin{landscape}
\begin{table}[h]
    \caption{Classification results for LP vs. SDP comparison in terms of accuracy. Bold entries indicate the best feature-based accuracy among XGB, GBC and RFC within each feasible setup. GNN-C is reported separately as a setup-independent benchmark.}
    \label{tab:classification-LPvsSDP}
    \hspace{-3cm} 
\begin{tabular}{|c|c|c|c|c|ccc|ccc|ccc|}
\hline
  &  & Training Set & Test Set 
  & \multicolumn{1}{c|}{GNN-C}
  & \multicolumn{3}{c|}{XGB-C}
  & \multicolumn{3}{c|}{GBC}
  & \multicolumn{3}{c|}{RFC} \\
\hline
Type & ID & $n,m,N$ & $n,m,N$ & 
   & fDD & sDD & DI 
   & fDD & sDD & DI
   & fDD & sDD & DI \\
\hline

\multirow{12}{*}{\rotatebox[origin=c]{90}{Synthetic}}
& 1  & 5,2,40K & 5,2,10K & 0.9485 & 0.9356 & 0.9348 & 0.9340 & 0.9378 & \textbf{0.9377} & \textbf{0.9380} & \textbf{0.9383} & 0.9322 & 0.9357 \\ \cline{2-14}

& 2  & 5,3,40K & 5,3,10K & 0.9380 & 0.9322 & 0.9307 & 0.9321 & 0.9315 & \textbf{0.9331} & 0.9334 & \textbf{0.9324} & 0.9307 & \textbf{0.9339} \\ \cline{2-14}

& 3  & 5,101,400 & 5,101,100 & 0.9048 & \textbf{0.8646} & \textbf{0.8542} & \textbf{0.8958} & 0.8437 & \textbf{0.8542} & 0.8750 & 0.854 & 0.823 & NA \\ \cline{2-14}

& 4  & 10,2,40K & 10,2,10K & 0.9800 & 0.9820 & 0.9823 & 0.9817 & 0.9826 & \textbf{0.9830} & 0.9823 & \textbf{0.9830} & 0.9823 & \textbf{0.9824} \\ \cline{2-14}

& 5  & 10,3,40K & 10,3,10K & 0.9693 & 0.9816 & 0.9811 & \textbf{0.9826} & \textbf{0.9821} & 0.9820 & 0.9821 & 0.9819 & \textbf{0.9822} & 0.9822 \\ \cline{2-14}

& 6  & 5,3,40K & 5,11,10K & 0.9479 & NA & NA & \textbf{0.9217} & NA & NA & 0.9184 & NA & NA & 0.9204 \\ \cline{2-14}

& 7  & 5,3,40K + 10,3,40K & 5,3,10K + 10,3,10K & 0.6856 & NA & 0.9754 & 0.9787 & NA & 0.9595 & 0.9600 & NA & \textbf{0.9996} & \textbf{0.9945}  \\ \cline{2-14}

& 8  & 5,3,50K + 10,3,50K & 20,3,10K & 0.9816 & NA & 0.9912 & 0.9905 & NA & \textbf{0.9953} & \textbf{0.9943} & NA & 0.9945 & 0.9932 \\ \cline{2-14}

& 9  & 10,3,50K & 5,11,10K & 0.9029 & NA & NA & \textbf{0.9042} & NA & NA & 0.9004 & NA & NA & 0.9009 \\ \cline{2-14}

& 10 & 5,11,10K & 10,3,50K & 0.9217 & NA & NA & 0.9069 & NA & NA & 0.9217 & NA & NA & \textbf{0.9412} \\ \cline{2-14}

& 11 & 10,3,50K & 5,101,500 & 0.8420 & NA & NA & \textbf{0.8680} & NA & NA & 0.8640 & NA & NA & \textbf{0.8680} \\ \cline{2-14}

\hline
\multirow{6}{*}{\rotatebox[origin=c]{90}{Synthetic$'$}}
& 12 & 5,2,40K & 5,2,10K & 0.9863 & \textbf{0.9855} & 0.9840 & 0.9828 & 0.9853 & \textbf{0.9845} & \textbf{0.9836} & 0.9854 & 0.9821 & 0.9828 \\ \cline{2-14}

& 13 & 5,101,400 & 5,101,100 & 0.9524 & 0.9380 & \textbf{0.9690} & 0.9270 & \textbf{0.9580} & \textbf{0.9690} & 0.9380 & 0.854 & 0.854 & \textbf{0.9690}  \\ \cline{2-14}

& 14 & 10,3,40K & 10,3,10K & 0.9803 & 0.9735 & 0.9719 & 0.9688 & 0.9737 & \textbf{0.9728} & \textbf{0.9709} & \textbf{0.9741} & 0.9693 & 0.9698 \\ \cline{2-14}

& 15 & 5,3,50K & 5,11,10K & 0.9031 & NA & NA & 0.9190 & NA & NA & 0.9290 & NA & NA & \textbf{0.9330} \\ \cline{2-14}

& 16 & 5,3,50K + 10,3,50K & 20,3,10K & 0.8945 &  NA & 0.9700 & 0.9600 & NA & \textbf{0.9710} & \textbf{0.9640} & NA & 0.968 & 0.9630 \\ \cline{2-14}

& 17 & 5,11,10K & 10,3,50K & 0.9601 & NA & NA & 0.9530 & NA & NA & 0.9620 & NA & NA & \textbf{0.9680} \\ \cline{2-14}

\hline
\end{tabular}  
\end{table}

\end{landscape}

We compare the GNN-C against all feature-based classifiers for each experiment.  
GNN-C achieves the highest overall accuracy in 6 of the 17 experiments (IDs 1, 2, 3, 6, 12 and 14), which shows that the graph-based approach can be effective in selected settings. 
 However, the performance of GNN-C is not consistent across the full experiments. For example, experiment 7 is one of the most structurally heterogeneous cases, and the GNN-C accuracy drops to 0.6856. There is also a noticeable gap for GNN between experiments 15 and 16. Therefore, we can conclude that GNN-C is more sensitive to increases in structural complexity and variation. From a practical standpoint, training GNN-C is also much slower: in our implementation, the training is much slower than the DI-based classifiers, and this gap tends to increase with problem complexity. Therefore, although GNN-C performs well in some experiments, the feature-based classifiers, especially under the DI setup, provide a more reliable balance of accuracy, applicability, and training time, and this pattern is also recognizable in other tables.

\paragraph{Regression Model Performance}

We use a custom-defined performance metric, called \textit{r-accuracy}, in the regression model, which is computed as follows:
Let  a vector of true labels \( y \) and a vector of predicted labels \( \hat{y} \) be given. Then,  we define the performance metric {r-accuracy} as 
\begin{equation*}
\text{r-accuracy}(y, \hat{y}) = \frac{1}{N} \sum_{\iota=1}^{N} \mathbf{1}\left\{ \text{sign}(y_\iota) = \text{sign}(\hat{y}_\iota) \, \text{or} \, \left(-\epsilon< y_\iota < 0 \, \text{and} \, \hat{y}_\iota > 0 \right) \right\},
\end{equation*}where  \( \epsilon \) is a predefined threshold for near-zero values. Our purpose in developing this performance metric is to consider predictions with the same sign as the actual relative difference as correctly categorized, regardless of their magnitude, and to also accept LP predictions in cases where SDP is tighter than LP by at most $\epsilon$ (the logic here is that  one might opt for a slightly weaker LP relaxation since it is computationally less demanding).  This metric also provides the advantage of making rough comparisons with classification accuracy. In this study, $\epsilon$ is chosen as $10^{-2}$.

We report  {r-accuracy}  as the main performance criterion in the Regression experiments reported in Table~\ref{tab:regression-LPvsSDP} for these regressors: GNN-R, XGB-R, GBR and RFR. Our preliminary experiments have shown that Support Vector Regressor is again   consistently outperformed by these four classifiers.

For each feasible experiment-setup pair, we highlight in bold the best feature-based regressor among XGB-R, GBR and RFR. In this comparison, RFR attains the largest number of best values with 18 outright wins. However, the three regressors are close to each other, with averages of 0.9376, 0.9335 and 0.9342 for XGB-R, GBR and RFR, respectively. Hence, although RFR appears slightly stronger in terms of pairwise wins, all three feature-based regressors perform at a similarly high level overall.

\begin{landscape}
\begin{table}[H]
    \caption{Regression results for LP vs. SDP comparison in terms of r-accuracy. Bold entries indicate the best feature-based result among XGB-R, GBR, and RFR within each feasible setup. GNN-R is reported separately as a setup-independent benchmark.}
    \label{tab:regression-LPvsSDP}
    \hspace{-3cm}
\begin{tabular}{|c|c|c|c|c|ccc|ccc|ccc|}
\hline
& & Training Set & Test Set & \multicolumn{1}{c|}{GNN-R} & \multicolumn{3}{c|}{XGB-R} & \multicolumn{3}{c|}{GBR} & \multicolumn{3}{c|}{RFR} \\
\hline
Type & ID & $n,m,N$ & $n,m,N$ & 
   & fDD & sDD & DI 
   & fDD & sDD & DI
   & fDD & sDD & DI \\
\hline

\multirow{12}{*}{\rotatebox[origin=c]{90}{Synthetic}}
& 1  & 5,2,40K & 5,2,10K & 0.8762 & 0.9364 & 0.9351 & 0.9358 & 0.9347 & 0.9353 & 0.9341 & \textbf{0.9372} & \textbf{0.9358} & \textbf{0.9366} \\ \cline{2-14}

& 2  & 5,3,40K & 5,3,10K & 0.8858 & \textbf{0.9303} & \textbf{0.9297} & \textbf{0.9312} & 0.9254 & 0.9277 & 0.8870 & 0.9294 & 0.9295 & 0.9305 \\ \cline{2-14}

& 3  & 5,101,400 & 5,101,100 & 0.8000 & \textbf{0.8854} & \textbf{0.8854} & 0.8958 & 0.8646 & \textbf{0.8854} & \textbf{0.9271} & \textbf{0.8854} & \textbf{0.8854} & 0.9063 \\ \cline{2-14}

& 4  & 10,2,40K & 10,2,10K & 0.9607 & \textbf{0.9825} & 0.9825 & 0.9944 & 0.9817 & 0.9818 & \textbf{0.9945} & 0.9823 & \textbf{0.9826} & 0.9830 \\ \cline{2-14}

& 5  & 10,3,40K & 10,3,10K & 0.9335 & 0.9811 & 0.9811 & \textbf{0.9826} & 0.9784 & 0.9784 & 0.9810 & \textbf{0.9815} & \textbf{0.9815} & 0.9816 \\ \cline{2-14}

& 6  & 5,3,40K & 5,11,10K & 0.8667 & NA & NA & 0.9086 & NA & NA & 0.9115 & NA & NA & \textbf{0.9117} \\ \cline{2-14}

& 7  & 5,3,40K + 10,3,40K & 5,3,10K + 10,3,10K & 0.5967 & NA & 0.9646 & 0.9640 & NA & 0.9537 & 0.9565 & NA & \textbf{0.9705} & \textbf{0.9713} \\ \cline{2-14}

& 8  & 5,3,50K + 10,3,50K & 20,3,10K & 0.8634 & NA & 0.9929 & 0.9851 & NA & 0.9892 & 0.9897 & NA & \textbf{0.9946} & \textbf{0.9911} \\ \cline{2-14}

& 9  & 10,3,50K & 5,11,10K & 0.8887 & NA & NA & 0.9053 & NA & NA & \textbf{0.9064} & NA & NA & 0.9044 \\ \cline{2-14}

& 10 & 5,11,10K & 10,3,50K & 0.9001 & NA & NA & \textbf{0.9341} & NA & NA & 0.9232 & NA & NA & 0.9321 \\ \cline{2-14}

& 11 & 10,3,50K & 5,101,500 & 0.6400 & NA & NA & \textbf{0.8660} & NA & NA & \textbf{0.8660} & NA & NA & \textbf{0.8660} \\ \cline{2-14}

\hline
\multirow{6}{*}{\rotatebox[origin=c]{90}{Synthetic$'$}}
& 12 & 5,2,40K & 5,2,10K & 0.9666 & 0.9656 & 0.9676 & 0.9555 & 0.9390 & 0.9315 & 0.9358 & \textbf{0.9765} & \textbf{0.9755} & \textbf{0.9693} \\ \cline{2-14}
& 13 & 5,101,400 & 5,101,100 & 0.7524 & \textbf{0.9380} & 0.8650 & 0.8230 & 0.9230 & 0.8440 & \textbf{0.9270} & 0.917 & \textbf{0.906} & 0.8850 \\ \cline{2-14}
& 14 & 10,3,40K & 10,3,10K & 0.9367 & 0.9564 & 0.9529 & 0.9370 & 0.9317 & 0.9368 & 0.9270 & \textbf{0.9667} & \textbf{0.9656} & \textbf{0.9507} \\ \cline{2-14}

& 15  & 5,3,50K & 5,11,10K & 0.7313 & NA & NA & \textbf{0.9060} & NA & NA & 0.9000 & NA & NA & 0.7340 \\ \cline{2-14}

& 16  & 5,3,50K + 10,3,50K & 20,3,10K & 0.7555 & NA & \textbf{0.9490} & 0.9170 & NA & 0.9310 & \textbf{0.9300} & NA & 0.9450 & 0.8920 \\ \cline{2-14}

& 17 & 5,11,10K & 10,3,50K & 0.9040 & NA & NA & 0.9310 & NA & NA & \textbf{0.9370} & NA & NA & 0.8360 \\ \cline{2-14}

\hline
\end{tabular}  
\end{table}
\end{landscape}

Although the accuracy metric in the classification task and the r-accuracy metric in the regression task are only partially comparable, we think that the latter might provide additional insights for relaxation selection purposes. This is due to the fact that the computational effort is somehow incorporated in the calculation of the r-accuracy metric, which is clearly a very important consideration.

Results of the regression task are generally strong. In 33 of the 36 feasible experiment setup pairs, at least one feature-based regressor achieves an r-accuracy of at least 0.90. Similar to the classification results, the DI setup remains attractive from the standpoint of applicability.

In contrast to the classification task, GNN-R does not achieve the best result in any of the 17 experiments. Its performance is weaker and decreases in more challenging settings.


\subsubsection{LP vs. SDP$'$ Comparison}
\label{subsec:comp_results_LPvsSDPprime}

\revise{We next present the results of our   relaxation
selection experiments on synthetically generated QCQP instances,
comparing LP and SDP$'$ relaxations.
Tables~\ref{tab:classification-LPvsSDP$'$} and
\ref{tab:regression-LPvsSDP$'$} report the classification
and regression results, respectively.}

\paragraph{Classification Model Performance}
In the classification task (Table~\ref{tab:classification-LPvsSDP$'$}), GBC attains the largest number of best feature-based results, with 19 wins and one tie, while XGB and RFC each obtain 8 wins. The overall classification performance remains strong: in 32 of the 36 feasible experiment setup pairs, at least one feature-based classifier achieves an accuracy of at least 0.90. The only exceptions are Experiment 3 under all three setups and Experiment 11 under the DI setup.

\paragraph{Regression Model Performance}

For the regression task (Table~\ref{tab:regression-LPvsSDP$'$}), the results are more decisive. RFR clearly dominates the feature-based comparison with 23 wins out of 36, while GBR and XGB obtain 9 and 4 wins, respectively. For both the classification and regression tasks, the DI setup remains competitive, but when multiple setups are applicable, the best results often come from the more disaggregated setups. While GNN-C results remain close to other ones, GNN-R performs much more poorly in the regression task. In particular, GNN-R falls below 0.60 in seven experiments and below 0.20 in Experiments 12, 14, and 16, indicating highly unstable behavior. 

Overall, compared to SDP, SDP$'$ relaxation is a much stronger competitor for LP relaxation and can give much tighter bounds than SDP by making good use of variable bounds. Consequently, the predictive task becomes more challenging, since the two relaxations often yield similarly strong performance. Nevertheless, the results indicate that the proposed approach remains successful in identifying the superior relaxation.

\begin{landscape}

\begin{table}[H]
    \caption{Classification results for LP vs. SDP$'$ comparison in terms of accuracy. Bold entries indicate the best feature-based result among XGB-C, GBC, and RFC within each feasible setup. GNN-C is reported separately as a setup-independent benchmark.}
    \label{tab:classification-LPvsSDP$'$}
    \hspace{-3cm} 
\begin{tabular}{|c|c|c|c|c|ccc|ccc|ccc|}
\hline
  &  & Training Set & Test Set 
  & \multicolumn{1}{c|}{GNN-C}
  & \multicolumn{3}{c|}{XGB-C}
  & \multicolumn{3}{c|}{GBC}
  & \multicolumn{3}{c|}{RFC} \\
\hline
Type & ID & $n,m,N$ & $n,m,N$ & 
   & fDD & sDD & DI 
   & fDD & sDD & DI
   & fDD & sDD & DI \\
\hline

\multirow{12}{*}{\rotatebox[origin=c]{90}{Synthetic}}
& 1  & 5,2,40K & 5,2,10K & 0.9479 & 0.9320 & 0.9345 & 0.9329 & \textbf{0.9378} & \textbf{0.9357} & \textbf{0.9346} & 0.9340 & 0.9305 & 0.9315 \\ \cline{2-14}

& 2  & 5,3,40K & 5,3,10K & 0.9527 & 0.9372 & 0.9369 & 0.9373 & \textbf{0.9386} & \textbf{0.9386} & \textbf{0.9386} & 0.9380 & 0.9357 & 0.9369 \\ \cline{2-14}

& 3  & 5,101,400 & 5,101,100 & 0.9143 & \textbf{0.8958} & \textbf{0.8854} & 0.8850 & 0.8854 & 0.8646 & 0.875 & 0.7396 & 0.8229 & \textbf{0.8960} \\ \cline{2-14}

& 4  & 10,2,40K & 10,2,10K & 0.9661 & \textbf{0.9678} & \textbf{0.9663} & 0.9664 & 0.9666 & \textbf{0.9663} & \textbf{0.9679} & 0.9660 & 0.9655 & 0.9641 \\ \cline{2-14}

& 5  & 10,3,40K & 10,3,10K & 0.9753 & 0.9737 & 0.9742 & 0.9746 & \textbf{0.9752} & \textbf{0.9762} & \textbf{0.9756} & 0.9742 & 0.9739 & 0.9740 \\ \cline{2-14}

& 6  & 5,3,40K & 5,11,10K & 0.8695 & NA & NA & 0.9439 & NA & NA & 0.8079 & NA & NA & \textbf{0.9463} \\ \cline{2-14}

& 7  & 5,3,40K + 10,3,40K & 5,3,10K + 10,3,10K & 0.7310 & NA & 0.9498 & \textbf{0.9473} & NA & \textbf{0.9645} & 0.9455 & NA & 0.9489 & 0.9454 \\ \cline{2-14}

& 8  & 5,3,50K + 10,3,50K & 20,3,10K & 0.9455 & NA & 0.9685 & 0.9728 & NA & 0.9595 & 0.9588 & NA & \textbf{0.9992} & \textbf{0.9984} \\ \cline{2-14}

& 9  & 10,3,50K & 5,11,10K & 0.9507 & NA & NA & 0.9435 & NA & NA & 0.8278 & NA & NA & \textbf{0.9447} \\ \cline{2-14}

& 10 & 5,11,10K & 10,3,50K & 0.9608 & NA & NA & \textbf{0.9639} & NA & NA & 0.9630 & NA & NA & 0.9599 \\ \cline{2-14}

& 11 & 10,3,50K & 5,101,500 & 0.8880 & NA & NA & \textbf{0.8820} & NA & NA & 0.7680 & NA & NA & 0.8780 \\ \cline{2-14}

\hline
\multirow{6}{*}{\rotatebox[origin=c]{90}{Synthetic$'$}}
& 12 & 5,2,40K & 5,2,10K & 0.9531 & 0.9798 & 0.9795 & \textbf{0.9797} & \textbf{0.9818} & \textbf{0.9816} & 0.9796 & 0.9803 & 0.9781 & 0.9786 \\ \cline{2-14}

& 13 & 5,101,400 & 5,101,100 & 0.9238 & 0.9583 & \textbf{0.9583} & 0.9375 & \textbf{0.9588} & 0.9480 & 0.9271 & 0.8438 & 0.8438 & \textbf{0.9583} \\ \cline{2-14}

& 14 & 10,3,40K & 10,3,10K & 0.9566 & 0.9677 & 0.9652  & 0.9628 & \textbf{0.9701} & \textbf{0.9659} & \textbf{0.9638} & 0.9685 & 0.9642 & 0.9629 \\ \cline{2-14}

& 15  & 5,3,50K & 5,11,10K & 0.9076 & NA & NA & 0.9069 & NA & NA & 0.9151 & NA & NA & \textbf{0.9181} \\ \cline{2-14}

& 16  & 5,3,50K + 10,3,50K & 20,3,10K & 0.9455 & NA & 0.9624 & 0.9494 & NA & \textbf{0.9636} & \textbf{0.9519} & NA & 0.9568 & 0.9509 \\ \cline{2-14}

& 17 & 5,11,10K & 10,3,50K & 0.9520 & NA & NA & 0.9389 & NA & NA & 0.9544 & NA & NA & \textbf{0.9607} \\ \cline{2-14}
\hline
\end{tabular}  
\end{table}

\end{landscape}

\begin{landscape}

\begin{table}[H]
    \caption{Regression results for LP vs. SDP$'$ comparison in terms of r-accuracy. Bold entries indicate the best feature-based result among XGB-R, GBR, and RFR within each feasible setup. GNN-R is reported separately as a setup-independent benchmark.}
    \label{tab:regression-LPvsSDP$'$}
    \hspace{-3cm}
\begin{tabular}{|c|c|c|c|c|ccc|ccc|ccc|}
\hline
& & Training Set & Test Set & \multicolumn{1}{c|}{GNN-R} & \multicolumn{3}{c|}{XGB-R} & \multicolumn{3}{c|}{GBR} & \multicolumn{3}{c|}{RFR} \\
\hline
Type & ID & $n,m,N$ & $n,m,N$ & 
   & fDD & sDD & DI 
   & fDD & sDD & DI
   & fDD & sDD & DI \\
\hline

\multirow{12}{*}{\rotatebox[origin=c]{90}{Synthetic}}
& 1  & 5,2,40K & 5,2,10K & 0.7954 & 0.9308 & 0.9305 & 0.9289 & \textbf{0.9323} & \textbf{0.9311} & 0.9258 & 0.9312 & 0.9303 & \textbf{0.9308} \\ \cline{2-14}

& 2  & 5,3,40K & 5,3,10K & 0.9306 & 0.9346 & 0.9354 & \textbf{0.9363} & 0.9352 & 0.9332 & 0.9350 & \textbf{0.9378} & \textbf{0.9363} & 0.9360 \\ \cline{2-14}

& 3  & 5,101,400 & 5,101,100 & 0.8190 & 0.7708 & 0.7708 & \textbf{0.8958} & 0.8750 & 0.8750 & 0.8750 & \textbf{0.8854} & \textbf{0.8958} & 0.8854 \\ \cline{2-14}

& 4  & 10,2,40K & 10,2,10K & 0.9030 & 0.9628 & 0.9633 & 0.9610 & \textbf{0.9805} & 0.9573 & 0.9590 & 0.9619 & \textbf{0.9636} & \textbf{0.9625} \\ \cline{2-14}

& 5  & 10,3,40K & 10,3,10K & 0.9514 & \textbf{0.9730} & \textbf{0.9726} & 0.9707 & 0.9682 & 0.9672 & 0.9677 & 0.9728 & 0.9717 & \textbf{0.9724} \\ \cline{2-14}

& 6  & 5,3,40K & 5,11,10K & 0.5019 & NA & NA & 0.9321 & NA & NA & 0.9426 & NA & NA & \textbf{0.9503} \\ \cline{2-14}

& 7  & 5,3,40K + 10,3,40K & 5,3,10K + 10,3,10K & 0.4935 & NA & 0.9530 & 0.7885 & NA & 0.9485 & 0.6620 & NA & \textbf{0.9605} & \textbf{0.9456} \\ \cline{2-14}

& 8  & 5,3,50K + 10,3,50K & 20,3,10K & 0.5570 & NA & 0.9568 & 0.9562 & NA & 0.9525 & 0.9488 & NA & \textbf{0.9694} & \textbf{0.9674} \\ \cline{2-14}

& 9  & 10,3,50K & 5,11,10K & 0.9472 & NA & NA & 0.9104 & NA & NA & \textbf{0.9430} & NA & NA & 0.9408 \\ \cline{2-14}

& 10 & 5,11,10K & 10,3,50K & 0.8140 & NA & NA & 0.9561 & NA & NA & \textbf{0.9613} & NA & NA & 0.9580 \\ \cline{2-14}

& 11 & 10,3,50K & 5,101,500 & 0.4400 & NA & NA & 0.8260 & NA & NA & \textbf{0.8780} & NA & NA & 0.8740 \\ \cline{2-14}

\hline
\multirow{6}{*}{\rotatebox[origin=c]{90}{Synthetic$'$}}
& 12 & 5,2,40K & 5,2,10K & 0.1614 & 0.9598 & 0.9626 & 0.9416 & 0.9604 & 0.9608 & 0.9445 & \textbf{0.9679} & \textbf{0.9688} & \textbf{0.9533} \\ \cline{2-14}

& 13 & 5,101,400 & 5,101,100 & 0.6000 & 0.7396 & 0.7500 & 0.8230 & 0.8854 & \textbf{0.9271} & 0.8438 & \textbf{0.9167} & 0.9063 & \textbf{0.8542} \\ \cline{2-14}

& 14 & 10,3,40K & 10,3,10K & 0.1858 & 0.9481 & 0.9443 & 0.9153 & 0.9409 & 0.9395 & 0.9092 & \textbf{0.9527} & \textbf{0.9448} & \textbf{0.9230} \\ \cline{2-14}

& 15  & 5,3,50K & 5,11,10K & 0.9164 & NA & NA & 0.7459 & NA & NA & \textbf{0.8285} & NA & NA & 0.7339 \\ \cline{2-14}

& 16  & 5,3,50K + 10,3,50K & 20,3,10K & 0.1895 & NA & 0.9377 & 0.8826 & NA & 0.9370 & \textbf{0.8951} & NA & \textbf{0.9446} & 0.8923 \\ \cline{2-14}

& 17 & 5,11,10K & 10,3,50K & 0.7360 & NA & NA & 0.8084 & NA & NA & 0.8892 & NA & NA & \textbf{0.9229} \\ \cline{2-14}

\hline
\end{tabular}  
\end{table}
\end{landscape}

\subsection{Benchmark Instances}
\label{subsec:benchmark_instances}

\revise{We now evaluate our relaxation selection framework on benchmark instances from MINLPLib and QPLIB. Since these instances vary in the number of variables and constraints, we use  DI  features for the feature-based models and the static tripartite graph representation for GNNs. We present relaxation selection results for LP vs.\ SDP and LP vs.\ SDP$'$, followed by relaxation selection results within spatial B\&B   for LP vs.\ SDP$'$.}

\revise{We consider two training approaches: Combined-to-Combined, where models are trained and tested on different subsets of the combined MINLPLib and QPLIB dataset, and Synthetic$'$-to-Combined, where models are trained on synthetic instances and evaluated on the combined benchmarks. This allows us to examine both benchmark-trained performance and transfer from synthetic to benchmark instances.}

\subsubsection{Relaxation Selection without Branch-and-Bound}

\revise{We first investigate relaxation selection between LP and SDP, and between LP and SDP$'$ without B\&B. Tables~\ref{tab:classification-MINLPLib} and~\ref{tab:regression-MINLPLib} report the classification and regression results using targets constructed from the corresponding root node lower bounds. Alongside the feature-based models, we evaluate GNN-C and GNN-R using the static tripartite graph representation.}

\revise{For classification, feature-based models trained on the combined benchmarks achieve strong performance in both comparisons. XGB-C obtains the highest accuracies, 0.9604 for LP vs.\ SDP and 0.9406 for LP vs.\ SDP$'$. RFC and GBC also exceed 0.89 in both comparisons, while GNN-C achieves 0.9208 and 0.8713, respectively. Thus, the feature-based models outperform GNN-C in this setting, indicating that DI features are useful for identifying the stronger root node relaxation.}

\revise{The differences become more pronounced in the Synthetic$'$-to-Combined experiments. For LP vs.\ SDP, XGB-C achieves the highest accuracy of 0.8165; for LP vs.\ SDP$'$, XGB-C and GBC reach 0.9209 and 0.9173, respectively. In contrast, GNN-C achieves only 0.4820 and 0.4460. Thus, under the current training setup, the DI feature representation transfers more effectively from synthetic to benchmark instances than the learned graph representation.}

\begin{table}[H]
    \caption{Root node classification results on the combined MINLPLib and QPLIB benchmark dataset for LP vs.\ SDP and LP vs.\ SDP$'$ comparisons in terms of accuracy. Bold entries indicate the best result across all models.}
    \label{tab:classification-MINLPLib}
    \centering
    \resizebox{\textwidth}{!}{
    \begin{tabular}{|m{2.5cm}|m{2.8cm}|c|c|c|c|}
    \hline
    \textbf{Type} & \textbf{Training Set}
    & \textbf{GNN-C} & \textbf{XGB-C} & \textbf{GBC} & \textbf{RFC} \\
    \hline
    LP vs.\ SDP & Combined
    & 0.9208 & \textbf{0.9604} & 0.9406 & 0.9505 \\
    \hline
    LP vs.\ SDP & Synthetic$'$
    & 0.4820 & \textbf{0.8165} & 0.7374 & 0.7986 \\
    \hline
    LP vs.\ SDP$'$ & Combined
    & 0.8713 & \textbf{0.9406} & 0.9010 & 0.8911 \\
    \hline
    LP vs.\ SDP$'$ & Synthetic$'$
    & 0.4460 & \textbf{0.9209} & 0.9173 & 0.8813 \\
    \hline
    \end{tabular}
    }
\end{table}

\revise{Table~\ref{tab:regression-MINLPLib} presents the regression results. In the Combined-to-Combined setting, RFR achieves the highest r-accuracy of 0.9307 for LP vs.\ SDP, while GBR leads LP vs.\ SDP$'$ with 0.9109. GNN-R achieves 0.8911 and 0.8515, respectively, comparable to some feature-based models but below the best results.}
\revise{For Synthetic$'$-to-Combined regression, RFR obtains the highest r-accuracy in both comparisons, reaching 0.8885 for LP vs.\ SDP and 0.7230 for LP vs.\ SDP$'$. GNN-R achieves 0.6655 and 0.5863, remaining below all feature-based regressors. Its predictions also exhibit very large absolute errors, indicating difficulty in estimating the magnitude of the relative bound difference.}

\revise{Interestingly, when trained on Synthetic$'$, feature-based classifiers perform better on LP vs.\ SDP$'$ than on LP vs.\ SDP, whereas regression models show the opposite pattern. Classification identifies the favorable relaxation, while regression also estimates the magnitude of the relative bound difference; these tasks need not be equally difficult across comparisons.}

\begin{table}[H]
    \caption{Root node regression results on the combined MINLPLib and QPLIB benchmark dataset for LP vs.\ SDP and LP vs.\ SDP$'$ comparisons in terms of r-accuracy. Bold entries indicate the best result across all models.}
    \label{tab:regression-MINLPLib}
    \centering
    \resizebox{\textwidth}{!}{
    \begin{tabular}{|m{2.5cm}|m{2.8cm}|c|c|c|c|}
    \hline
    \textbf{Type} & \textbf{Training Set}
    & \textbf{GNN-R} & \textbf{XGB-R} & \textbf{GBR} & \textbf{RFR} \\
    \hline
    LP vs.\ SDP & Combined
    & 0.8911 & 0.8911 & 0.9109 & \textbf{0.9307} \\
    \hline
    LP vs.\ SDP & Synthetic$'$
    & 0.6655 & 0.8597 & 0.7806 & \textbf{0.8885} \\
    \hline
    LP vs.\ SDP$'$ & Combined
    & 0.8515 & 0.8515 & \textbf{0.9109} & 0.8812 \\
    \hline
    LP vs.\ SDP$'$ & Synthetic$'$
    & 0.5863 & 0.6547 & 0.6655 & \textbf{0.7230} \\
    \hline
    \end{tabular}
    }
\end{table}

\revise{Feature-based models perform well in  several Synthetic$'$-to-Combined settings. GNNs are competitive when trained on combined benchmarks but generalize less effectively from synthetic instances in both classification and regression. To further investigate the GNN results, we examined instance-level predictions on the benchmark dataset. In Synthetic$'$-to-Combined classification, GNN-C predicts LP for only 35.3\% and 31.7\% of instances in the LP--SDP and LP--SDP$'$ comparisons, respectively, whereas the corresponding true proportions are 80.6\% and 64.7\%. Both accuracies fall below the majority-class baselines. GNN-R also exhibits weaker cross-dataset performance, with r-accuracy decreasing from 0.8911 to 0.6655 for LP--SDP and from 0.8515 to 0.5863 for LP--SDP$'$. These findings suggest limited transfer of the learned GNN representations from synthetic to benchmark instances. Additional diagnostics are provided in Appendix~\ref{app:gnn_diagnostics}.}

\subsubsection{Relaxation Selection within Spatial Branch-and-Bound }

\revise{We now extend the evaluation to spatial B\&B, comparing global lower bounds from independently generated LP and SDP$'$ trees under the time-calibrated budgets introduced earlier. Unlike root node targets, these measure lower bound progress during the search. The combined benchmark contains 280 MINLPLib and QPLIB instances: 177 for training and 103 for testing in Combined-to-Combined experiments. Synthetic$'$-to-Combined models are trained on synthetic instances and evaluated on the combined benchmark.}

\revise{Tables~\ref{tab:classification-combined-bb} and~\ref{tab:regression-combined-bb} report the spatial B\&B classification and regression results. Feature-based models use 64 static DI features, while GNNs use the QCQP tripartite graph. Both representations are available before the search, so neither approach relies on observations collected during B\&B.}

\revise{For Combined-to-Combined classification, GBC and RFC achieve the highest accuracy of 0.9029, followed by XGB-C (0.8738) and GNN-C (0.8252). Thus, DI features remain informative when the target is based on global lower bound progress rather than the root node relaxation values.}

\revise{In Synthetic$'$-to-Combined classification, XGB-C achieves the highest accuracy of 0.9500, followed by GBC (0.9214) and RFC (0.8857). XGB-C and GBC improve over their Combined-to-Combined results, suggesting that Synthetic$'$ provides useful training information for spatial B\&B relaxation selection. In contrast, GNN-C drops to 0.2286, indicating that cross-dataset transfer depends strongly on the representation and learning architecture.}

\begin{table}[H]
    \caption{Spatial B\&B classification results on the combined MINLPLib and QPLIB benchmark dataset for LP vs.\ SDP$'$ relaxation selection in terms of accuracy. Bold entries indicate the best feature-based result among XGB-C, GBC and RFC.}
    \label{tab:classification-combined-bb}
    \centering
    \begin{tabular}{|m{3cm}|c|c|c|c|}
    \hline
    \textbf{Training Set}
    & \textbf{GNN-C}
    & \textbf{XGB-C}
    & \textbf{GBC}
    & \textbf{RFC} \\
    \hline
    Combined
    & 0.8252
    & 0.8738
    & \textbf{0.9029}
    & \textbf{0.9029} \\
    \hline
    Synthetic$'$
    & 0.2286
    & \textbf{0.9500}
    & 0.9214
    & 0.8857 \\
    \hline
    \end{tabular}
\end{table}

\revise{For Combined-to-Combined B\&B regression, GBR achieves the highest r-accuracy of 0.8544, followed closely by RFR (0.8447) and XGB-R (0.8350), while GNN-R achieves 0.7670. The three feature-based regressors perform similarly, with GBR providing the strongest result.}

\revise{With Synthetic$'$ training, B\&B regression performance decreases: RFR achieves the highest r-accuracy of 0.6893, followed by GBR (0.6393), XGB-R (0.5571), and GNN-R (0.5607). Unlike classification, synthetic-trained regressors do not improve over benchmark-trained models. Thus, identifying the favorable relaxation transfers more successfully than estimating the magnitude of the global lower bound difference.}

\revise{The GNN results likewise show weaker Synthetic$'$-to-Combined transfer, especially for classification. These experiments do not establish that handcrafted features are universally preferable, but the DI representation performs better in the benchmark transfer settings considered here.}

\begin{table}[H]
    \caption{Spatial B\&B regression results on the combined MINLPLib and QPLIB benchmark dataset for LP vs.\ SDP$'$ relaxation selection in terms of r-accuracy. Bold entries indicate the best feature-based result among XGB-R, GBR and RFR.}
    \label{tab:regression-combined-bb}
    \centering
    \begin{tabular}{|m{3cm}|c|c|c|c|}
    \hline
    \textbf{Training Set}
    & \textbf{GNN-R}
    & \textbf{XGB-R}
    & \textbf{GBR}
    & \textbf{RFR} \\
    \hline
    Combined
    & 0.7670
    & 0.8350
    & \textbf{0.8544}
    & 0.8447 \\
    \hline
    Synthetic$'$
    & 0.5607
    & 0.5571
    & 0.6393
    & \textbf{0.6893} \\
    \hline
    \end{tabular}
\end{table}

\revise{Comparing root node and spatial B\&B classification for LP vs.\ SDP$'$, the best accuracies are 0.9406 and 0.9029 in the Combined-to-Combined setting, and 0.9209 and 0.9500 in the Synthetic$'$-to-Combined setting, respectively. The strong results in both settings suggest that   DI features capture information useful not only for initial relaxation strength but also for global lower bound progress under time-calibrated B\&B budgets.}

\revise{The regression models achieve lower r-accuracy values in the spatial B\&B experiments compared with the root node experiments. For LP vs.\ SDP$'$, the best r-accuracy decreases from 0.9109 at the root node to 0.8544 in B\&B in the Combined-to-Combined setting, and from 0.7230 to 0.6893 in the Synthetic$'$-to-Combined setting. This difference may be explained by the additional factors affecting relaxation performance throughout the B\&B search. At the root node, the target is determined by the difference between the initial relaxation bounds. In B\&B, the relative performance of LP and SDP$'$ also depends on the evolution of the global lower bounds as the variable domains are partitioned, as well as the computational effort required to solve the relaxations. Since the models rely only on static features extracted from the original QCQP, they do not directly capture how the two relaxations behave at subsequent nodes. This may limit their ability to predict the relative bound difference in B\&B as accurately as at the root node. Incorporating dynamic features collected during the early nodes of B\&B search, such as changes in global lower bounds and node-level relaxation performance, could therefore be a promising direction for future research. Such features may provide additional information about the progression of the search and potentially improve regression performance.}

\revise{Overall, learning-based relaxation selection remains promising in spatial B\&B.   DI features provide strong classification results, while predicting the magnitude of lower bound differences is more difficult; GNNs are weaker in the benchmark transfer experiments, particularly with synthetic training. 
}

\subsubsection{\revise{Gap--Time Trade-off Analysis}}

\revise{We next study the trade-off between relaxation quality and computational time for the LP and SDP$'$ relaxation-selection strategies. Figure~\ref{fig:gap-time-tradeoff} presents the gap--time trade-off for root node and spatial B\&B classification and regression. We focus on LP vs.\ SDP$'$, since SDP$'$ is the SDP-based relaxation used in our spatial B\&B experiments. All learned models are trained on Synthetic$'$ and evaluated on the combined MINLPLib and QPLIB benchmark instances. The horizontal and vertical axes report the average normalized optimality gap and average computational time, respectively; smaller values are preferred.}

\begin{figure}[H]
\centering
\begin{tikzpicture}[scale=0.95]
\pgfplotsset{
  idealstyle/.style={only marks,mark=asterisk,mark size=4pt,color=black,thick},
  bfstyle/.style={only marks,mark=x,mark size=4pt,color=orange!90!black,thick},
  xgbstyle/.style={only marks,mark=square,mark size=4pt,color=blue!80!black},
  gbstyle/.style={only marks,mark=o,mark size=4pt,color=green!60!black,mark options={fill=white},thick},
  rfstyle/.style={only marks,mark=diamond,mark size=4pt,color=red!80!black},
  gnnstyle/.style={only marks,mark=triangle*,mark size=4pt,color=violet!80!black},
  alpstyle/.style={only marks,mark=+,mark size=4pt,color=brown!80!black,thick},
  asdpstyle/.style={only marks,mark=star,mark size=4pt,color=gray!75!black}
}
\begin{groupplot}[
  group style={group size=2 by 2,horizontal sep=1.7cm,vertical sep=1.8cm},
  width=0.45\textwidth,height=0.45\textwidth,
  xmin=0.26,xmax=0.54,
  xlabel={Average normalized gap $\bar{g}$},
  ylabel={Average total time (s)},
  grid=both,
  grid style={line width=0.3pt,draw=black!10},
  major grid style={line width=0.5pt,draw=black!20},
  title style={font=\small},label style={font=\small},
  tick label style={font=\scriptsize},
  scaled x ticks=false,scaled y ticks=false,
  xtick={0.30,0.35,0.40,0.45,0.50},
  xticklabel style={/pgf/number format/fixed,/pgf/number format/precision=2},
  yticklabel style={/pgf/number format/fixed,/pgf/number format/precision=1},
  legend columns=4,
  legend style={font=\scriptsize,cells={anchor=west},at={(-0.15,-1.75)},anchor=north},
  legend image post style={scale=0.75}
]
\nextgroupplot[title={Root node -- classification},ymin=0,ymax=16]
\addplot[idealstyle] coordinates {(0.319300,0.076862)};
\addplot[bfstyle] coordinates {(0.319300,0.190133)};
\addplot[alpstyle] coordinates {(0.500422,0.104791)};
\addplot[asdpstyle] coordinates {(0.366413,0.085342)};
\addplot[xgbstyle] coordinates {(0.328506,9.020914)};
\addplot[gbstyle] coordinates {(0.324119,9.018625)};
\addplot[rfstyle] coordinates {(0.333897,9.014365)};
\addplot[gnnstyle] coordinates {(0.398098,0.915228)};

\nextgroupplot[title={Root node -- regression},ymin=0,ymax=16]
\addplot[idealstyle] coordinates {(0.319300,0.076862)};\addlegendentry{Ideal}
\addplot[bfstyle] coordinates {(0.319300,0.190133)};\addlegendentry{Brute Force}
\addplot[alpstyle] coordinates {(0.500422,0.104791)};\addlegendentry{Always LP}
\addplot[asdpstyle] coordinates {(0.366413,0.085342)};\addlegendentry{Always SDP$'$}
\addplot[xgbstyle] coordinates {(0.364513,9.029017)};\addlegendentry{XGB}
\addplot[gbstyle] coordinates {(0.355028,9.026646)};\addlegendentry{GB}
\addplot[rfstyle] coordinates {(0.344040,9.026465)};\addlegendentry{RF}
\addplot[gnnstyle] coordinates {(0.383441,0.936243)};\addlegendentry{GNN}

\nextgroupplot[title={Spatial B\&B -- classification},ymin=0,ymax=16]
\addplot[idealstyle] coordinates {(0.292861,4.613668)};
\addplot[bfstyle] coordinates {(0.292861,14.688360)};
\addplot[alpstyle] coordinates {(0.473353,10.531494)};
\addplot[asdpstyle] coordinates {(0.358075,4.156866)};
\addplot[xgbstyle] coordinates {(0.301755,13.631638)};
\addplot[gbstyle] coordinates {(0.304641,13.423633)};
\addplot[rfstyle] coordinates {(0.313384,13.712614)};
\addplot[gnnstyle] coordinates {(0.511259,10.733635)};

\nextgroupplot[title={Spatial B\&B -- regression},ymin=0,ymax=16]
\addplot[idealstyle] coordinates {(0.292861,4.613668)};
\addplot[bfstyle] coordinates {(0.292861,14.688360)};
\addplot[alpstyle] coordinates {(0.473353,10.531494)};
\addplot[asdpstyle] coordinates {(0.358075,4.156866)};
\addplot[xgbstyle] coordinates {(0.349897,13.025849)};
\addplot[gbstyle] coordinates {(0.335129,13.006129)};
\addplot[rfstyle] coordinates {(0.325403,13.002474)};
\addplot[gnnstyle] coordinates {(0.466259,11.119409)};
\end{groupplot}
\end{tikzpicture}
\caption{\revise{Gap--time trade-off for LP versus SDP$'$ relaxation selection
 on the combined MINLPLib and QPLIB benchmark. 
 GB denotes GBC
 and GBR, RF denotes RFC and RFR, and GNN denotes GNN-C and GNN-R in the
 corresponding panels.
 Lower values on both axes are preferable.}}
\label{fig:gap-time-tradeoff}
\end{figure}

\revise{For instance $\iota$, let $z_{\mathrm{UB}}^{(\iota)}$ denote the common feasible upper bound and $z_{\mathrm{sel}}^{(\iota)}$ the lower bound obtained by the selected relaxation. We calculate the normalized gap as}
\begin{equation}
 g^{(\iota)}=\frac{ z_{\mathrm{UB}}^{(\iota)}-z_{\mathrm{sel}}^{(\iota)} }
 {|z_{\mathrm{UB}}^{(\iota)}|+|z_{\mathrm{sel}}^{(\iota)}|+10^{-3}}.
 \label{eq:norm-gap}
\end{equation}
\revise{
}

\revise{The Always LP and Always SDP$'$ strategies solve the LP and SDP$'$ relaxations, respectively, for every instance. Their computational times correspond to the solution times of these relaxations.
 The Ideal strategy selects the stronger lower bound for each instance and counts only
 the selected solver time; it is an oracle benchmark, not a directly implementable
 strategy. Brute Force solves both relaxations and uses the sum of their solution
 times. For the feature-based models (XGB, GB and RF), we add the recorded time
 for computing the 64 DI features once per instance to the selected relaxation's
 solver time. The reported times do not separately include the feature-based
 model prediction time.}

\revise{For the GNN models, computational time includes graph construction, graph-feature normalization, model prediction, and the selected relaxation's solution time. Graph construction converts each QCQP into a tripartite graph of variable, quadratic-term, and constraint nodes; normalization uses statistics from the Synthetic$'$ training data. We measure each component separately on the CPU, with two warm-up runs followed by three timed runs per instance, using their median to reduce fluctuations. Model training and computation of normalization statistics are excluded. Thus, the reported times estimate per-instance cost in our experimental setup.
}

\revise{The average time of a strategy is the arithmetic mean of its per-instance total times. root node times correspond to the root relaxation solves, while B\&B times accumulate solution times over the selected LP or SDP$'$ tree. In spatial B\&B, each model selects a relaxation once per QCQP and uses it throughout the tree, so feature computation and prediction are counted only~once.}

\revise{At the root node, GBC obtains the smallest average gap among learned classifiers (0.3241), while RFR performs best among learned regressors (0.3440). The feature-based models generally provide stronger bounds than GNNs, with gaps close to Brute Force (0.3193). However, their computational times are approximately 9 seconds due to DI feature computation. GNN graph construction, normalization, and prediction require approximately 0.83--0.85 seconds, still more than the average Brute Force solver time (0.19 seconds). Thus, the learned models do not provide a computational advantage at the root node.}

\revise{The spatial B\&B experiments were conducted with a limited node budget because solving the relaxations repeatedly through the optimization libraries was computationally expensive. With this limited budget, the cost of feature extraction and prediction is relatively high compared to the time spent solving the relaxations. Using a larger node budget could allow better relaxation choices to reduce the number of explored nodes and improve computational time. However, this remains to be a future research.}

\revise{In spatial B\&B, XGB-C achieves the smallest average gap among learned classifiers (0.3018), while RFR performs best among learned regressors (0.3254). Both improve upon Always SDP$'$ (0.3581), and XGB-C approaches Brute Force and Ideal (0.2929). Feature-based models require approximately 13 seconds on average, less than Brute Force (14.69 seconds). 
}

\revise{Comparing the root node and spatial B\&B results reveals an important difference in the performance of the learned models. The GNN average gaps increase from 0.3981 to 0.5113 for classification and from 0.3834 to 0.4663 for regression, while the feature-based models maintain smaller gaps. One possible explanation is how these approaches represent QCQP instances. The DI features explicitly include global properties of the quadratic coefficient matrices, such as eigenvalue information, which helps describe their structure and curvature. In contrast, GNNs process individual variables, quadratic terms, constraints, and their connections through message passing. Although these components contain information about the original matrices, learning global matrix properties from them may be more difficult than using features that describe these properties directly. This may partly explain the better performance of feature-based models in both root node and B\&B experiments. 
}


\subsection{Feature Importance Analysis}
\label{subsec:feature_importance}

\revise{For the benchmark instances, we assess the importance of the 64 dimension-independent (DI) features using permutation importance: the decrease in evaluation performance when a feature is randomly permuted while the trained model remains fixed. Larger positive values indicate greater model reliance on that feature. We report selected features from representative experiments. 
}

\subsubsection{Feature Importance in Relaxation Selection  without Spatial B\&B}

\revise{Table~\ref{tab:root-feature-importance} shows that  spectral and sparsity characteristics of the objective function matrix $A_0$  are prominent across root node comparisons, although their importance varies by model and training dataset. BoundsExist is also informative for Synthetic$'$-to-Combined LP vs.\ SDP$'$ classification, highlighting the potential relevance of variable-bound availability in cross-dataset prediction.}

\begin{table}[H]
\caption{Selected permutation feature importance results for
root node relaxation selection. The reported values measure
the decrease in the corresponding evaluation score when
a feature is permuted.}
\label{tab:root-feature-importance}
\centering
\small
\begin{tabular}{|l|l|l|l|c|}
\hline
\textbf{Comparison} &
\textbf{Training} &
\textbf{Model} &
\textbf{Feature} &
\textbf{Importance} \\
\hline
LP vs.\ SDP & Combined & XGB-C
& $\eta_{A_0}/n$ & 0.1111 \\
& & & $\mathrm{PH}_{A_0}$ & 0.0676 \\
\hline
LP vs.\ SDP & Synthetic$'$ & RFR
& $\theta_{A_0}$ & 0.1302 \\
& & & $\widetilde{\lambda}_{\min}(A_0)$ & 0.0810 \\
\hline
LP vs.\ SDP$'$ & Combined & GBR
& $\mathrm{PH}_{A_0}$ & 0.2035 \\
& & & $\widetilde{\lambda}_{\min}(A_0)$ & 0.0919 \\
\hline
LP vs.\ SDP$'$ & Synthetic$'$ & XGB-C
& $\mathrm{PH}_{A_0}$ & 0.1579 \\
& & & BoundsExist & 0.0748 \\
\hline
\end{tabular}
\end{table}

\subsubsection{Feature Importance in Relaxation Selection  within Spatial B\&B}

\revise{Table~\ref{tab:bb-feature-importance} shows a similar pattern for spatial B\&B: objective spectral and sparsity features remain influential, while BoundsExist is among the two leading features in both Synthetic$'$-to-Combined classification and regression.}

\begin{table}[H]
\caption{Selected permutation feature importance results
for spatial B\&B relaxation selection between LP and SDP$'$.
The reported values measure the decrease in the
corresponding evaluation score when a feature is permuted.}
\label{tab:bb-feature-importance}
\centering
\small
\begin{tabular}{|l|l|l|c|}
\hline
\textbf{Experiment} &
\textbf{Model} &
\textbf{Feature} &
\textbf{Importance} \\
\hline
Combined / Classification & GBC
& $\widetilde{\lambda}_{\max}(A_0)$ & 0.0801 \\
& & $\eta_{A_0}/n$ & 0.0434 \\
\hline
Combined / Regression & RFR
& $\mathrm{PH}_{A_0}$ & 0.1112 \\
& & $\rho_{A_0}/n$ & 0.0462 \\
\hline
Synthetic$'$ / Classification & XGB-C
& $\mathrm{PH}_{A_0}$ & 0.1279 \\
& & BoundsExist & 0.1110 \\
\hline
Synthetic$'$ / Regression & RFR
& BoundsExist & 0.1198 \\
& & $\mathrm{PH}_{A_0}$ & 0.0307 \\
\hline
\end{tabular}
\end{table}

\subsection{Feature Ablation Study}
\label{subsec:feature_ablation}

\revise{For the benchmark instances, we complement permutation importance with ablation experiments that remove selected DI feature groups and retrain the models. We compare the full 64-feature representation against variants excluding objective-, constraint-, or bound-related features. Classification results use balanced accuracy to account for class imbalance, and regression results use the previously defined r-accuracy. The reported performance drop equals the full-feature score minus the ablated score.}

\subsubsection{Ablation Study for Relaxation Selection  without Spatial B\&B}

\revise{As Table~\ref{tab:root-ablation} shows, removing objective features reduces Combined-to-Combined XGB-C balanced accuracy from 0.9383 to 0.8755 and GBR r-accuracy from 0.9126 to 0.8155. In Synthetic$'$-to-Combined classification, excluding constraint features reduces GBC balanced accuracy from 0.9351 to 0.7611, whereas removing BoundsExist reduces XGB-C balanced accuracy from 0.9177 to 0.7736. Thus, the contributions of different feature groups depend on the learning model and training setting.}

\begin{table}[H]
\caption{Selected feature ablation results for root node
relaxation selection. 
}
\label{tab:root-ablation}
\centering
\small
\begin{tabular}{|l|l|l|c|c|c|}
\hline
\textbf{Setting} &
\textbf{Model} &
\textbf{Removed features} &
\textbf{Full} &
\textbf{Ablated} &
\textbf{Drop} \\
\hline
Combined / C & XGB-C
& Objective $A_0$ & 0.9383 & 0.8755 & 0.0628 \\
\hline
Combined / R & GBR
& Objective $A_0$ & 0.9126 & 0.8155 & 0.0971 \\
\hline
Synthetic$'$ / C & GBC
& Constraint features & 0.9351 & 0.7611 & 0.1740 \\
\hline
Synthetic$'$ / C & XGB-C
& BoundsExist & 0.9177 & 0.7736 & 0.1441 \\
\hline
\end{tabular}
\end{table}

\subsubsection{Ablation Study for Relaxation Selection  within Spatial B\&B}

\revise{Table~\ref{tab:bb-ablation} shows smaller losses from removing constraint or objective function features in Combined-to-Combined classification, but greater sensitivity in Synthetic$'$-to-Combined experiments. In particular, removing BoundsExist reduces XGB-C balanced accuracy by approximately 0.24, consistent with its high permutation importance. Removing objective features also decreases Synthetic$'$-to-Combined RFR r-accuracy from 0.6893 to 0.5179. Overall, the two analyses indicate that objective structure, constraints, and bound availability provide complementary predictive information, with their relative contributions varying across tasks and datasets.}

\begin{table}[H]
\caption{Selected feature ablation results for spatial
B\&B relaxation selection. 
}
\label{tab:bb-ablation}
\centering
\small
\begin{tabular}{|l|l|l|c|c|c|}
\hline
\textbf{Setting} &
\textbf{Model} &
\textbf{Removed features} &
\textbf{Full} &
\textbf{Ablated} &
\textbf{Drop} \\
\hline
Combined / C & GBC
& Constraint features & 0.8998 & 0.8590 & 0.0408 \\
\hline
Combined / C & GBC
& Objective shape & 0.8998 & 0.8581 & 0.0418 \\
\hline
Synthetic$'$ / C & XGB-C
& BoundsExist & 0.9430 & 0.7034 & 0.2397 \\
\hline
Synthetic$'$ / C & XGB-C
& Objective shape & 0.9430 & 0.7938 & 0.1492 \\
\hline
Synthetic$'$ / R & RFR
& Objective $A_0$ & 0.6893 & 0.5179 & 0.1714 \\
\hline
Synthetic$'$ / R & RFR
& Constraint features & 0.6893 & 0.5893 & 0.1000 \\
\hline
\end{tabular}
\end{table}

\section{Conclusions}
\label{sec:conclusions}

This paper focuses on analyzing and comparing the performance of three fundamental relaxations used in nonconvex QCQP problems, namely LP, SDP, and SDP$'$ relaxations. We formulate relaxation selection as a supervised learning task and develop both classification and regression models to predict the more favorable relaxation. Three distinct feature generation approaches successfully classify QCQP instances into SDP-favoring or LP-favoring categories. Our computational results show that the proposed feature-based framework can successfully predict the preferable relaxation without solving all alternatives in advance. In the classification tasks, Gradient Boosting is the strongest overall among the feature-based models, while in the regression tasks Random Forest is the strongest overall. Although the GNN classifier is competitive in some particular settings, its performance is less consistent across experiments, and the GNN regressor is generally weaker and more unstable. Therefore, the main practical conclusion of this study is that the feature-based framework provides a reliable approach for relaxation selection in nonconvex QCQPs.

The developed models provide a significant advantage by predicting the most beneficial relaxation type for a new QCQP instance without requiring users to test various relaxation methods. Another important outcome  is that the DI setup preserves strong predictive performance despite the aggregated features. This makes the approach attractive not only for synthetic experiments but also for benchmark and application-driven QCQP instances. Since the DI setup is independent of the instance size, it allows us to apply the same framework across QCQPs with different numbers of variables and constraints.

The proposed methodology also generalizes beyond synthetic data. In particular, our experiments on \revise{the combined MINLPLib and QPLIB benchmark} show that the proposed methodology transfers successfully to benchmark instances from several application domains, achieving around or above 90\% accuracy in the main classification settings. This supports the practical relevance of the framework beyond synthetic data and shows that the DI setup can be used on realistic QCQP instances arising from different problem families. \revise{We additionally evaluate relaxation selection within spatial B\&B, where the selected relaxation is used throughout the search tree.} We \revise{also} complement the prediction results with a gap-time trade-off analysis on \revise{the combined MINLPLib and QPLIB benchmark}, which incorporates solve times in addition to bound quality. This analysis shows that learned selectors can achieve a favorable balance between relaxation strength and computational cost while avoiding the expense of solving both relaxations. \revise{In the spatial B\&B experiments, feature-based selectors achieve lower bound quality close to that of solving both relaxations, although feature-computation cost and the limited node budget restrict the observed runtime savings.} Hence, the proposed framework is useful not only from a prediction standpoint but also from a computational efficiency perspective, \revise{specially for instances that may require a large spatial B\&B tree}.

There are several promising directions for future study. First, the relaxation runtimes could be incorporated more directly into the training objective so that prediction quality and computational efficiency are optimized jointly. Second, the practicality of the model could be improved by investigating the use of alternative features that are easier to compute, including eigenvalue estimates rather than precise calculations. Third, the accuracy of the framework could be improved by incorporating different sets for properties $\Pi$, statistics $Q$, and patterns $S$. Finally, \revise{the proposed relaxation selection framework can be evaluated on larger spatial B\&B trees with higher node budgets, where stronger relaxation choices may yield greater computational savings.}

\section*{Data Availability}

Our repository, which includes the codes to generate and solve instances \revise{as well as detailed statistics}, is available at 
\url{https://github.com/mugededeoglu/learning-relaxation}.

\revise{
\section*{Acknowledgments}
We sincerely thank the anonymous reviewers for their constructive comments, which helped us significantly improve our paper.
}

\section*{Statements and Declarations}

\subsection*{Competing interests}
The authors declare that they have no competing interests.

\appendix

\section{Instance Generation}
\label{sec:instance_generation}

In this section, we describe the process of generating QCQP instances. Our primary objective is to generate a diverse dataset to study how different properties of the instance structure affect the success of relaxations. We hypothesize that specific instance structures, such as the number of variables (\(n\)), the number of constraints (\(m\)) \cite{ghaddar2022learning,gonzalez2022polynomial}, matrix rank, and the existence of finite variable bounds, are important factors to consider. Additionally, the effect of the convex-concave nature of the problem has been discussed in previous works by \cite{fu1998approximation}, \cite{nemirovskii1999maximization}, and \cite{nesterov2000global}, and the presence of a single reverse convex constraint \cite{hillestad1980linear}.

The success of relaxations may be influenced by various factors, including  spectral properties, as discussed in \cite{lodi2022}. Other factors include the presence of bilinear constraints \cite{McCormick1976,torres1990linearization} and sparsity patterns (such as bipartite, tree, planar, chordal, and odd-cycle-rich graphs), as well as the diagonal nature of matrices \cite{kim2003exact,burer2018exact}.
 We aim to create a sufficiently diverse dataset so that we can systematically observe the effects of these instance structures. Section \ref{sec:sparsity} focuses on matrices that exhibit various sparsity patterns. We use the \texttt{Python NetworkX} package \cite{networkx} to create and manipulate graphs, including operations such as adding or removing nodes and edges. In Section \ref{sec:random}, we discuss the generation of \(n \times n\) random symmetric matrices with controlled properties. 
 In the final implementation, bilinear instances are generated by enforcing a zero diagonal on eligible support patterns, and mixed-size datasets are produced by combining instances with different values of \(n\) and~\(m\).
 
\begin{notation}
The symbol \(\unif(a, b)\) denotes a uniform distribution with bounds \(a\) and \(b\). 
\end{notation}

\subsection{Matrices with Different Sparsity Patterns}\label{sec:sparsity}

We now describe the generation of various types of $n \times n$ matrices used in our study. {These matrices include symmetric hollow, bipartite, tree, planar, chordal, and diagonal forms. }

\subsubsection{Bilinear (Zero-Diagonal) Matrices}

In the final implementation, bilinear matrices are generated by enforcing a zero diagonal on a weighted symmetric support pattern. We first generate a support matrix from an eligible sparsity family such as bipartite, tree, planar, chordal, odd-cycle-rich, symmetric, or eigenvalue-based structures. Then, random weights are assigned to the nonzero off-diagonal entries, and the diagonal entries are set exactly to zero. This construction yields bilinear-type quadratic forms while preserving the intended sparsity structure.

In some bounded bilinear instances, multiple matrices are generated from a shared support pattern in order to increase cross-constraint coupling. This makes the resulting instances more challenging for LP-based relaxations and helps diversify cases where SDP-based relaxations may perform better.

\subsubsection{Bipartite Graph Matrices}

We use the \texttt{bipartite} module of the \texttt{networkx} package to generate random bipartite graphs. First, we construct a bipartite graph \(G\) with \(k\) nodes in one partition and \(n-k\) nodes in the other partition. The adjacency matrix \(A\) of graph \(G\) is computed. Afterwards, \(A\) is converted into a weighted symmetric matrix \(D\), in which each edge (with a value of 1) is substituted with a random integer, indicating the distance between nodes. The latter two procedures are applicable to all subsequent graph types.

\begin{algorithm}[H]
\caption{Bipartite Matrix Generation.}\label{alg:bipartite_mtrx}
\begin{algorithmic}[1]
\Input $n$
\Output $D$: An $n \times n$ bipartite weighted symmetric matrix
\State Randomly choose $k$ such that $2 \leq k \leq n-2$
\State Randomly choose $e$ such that $2 \leq e \leq k \times (n - k)$
\State Generate a bipartite graph $G$ with $k$ and $n-k$ nodes and $e$ edges 
\State Compute the adjacency matrix $A$ of graph $G$ \label{step:adjacency} \Comment{Common Step 1}
\State Convert adjacency matrix $A$ into weighted symmetric matrix $D$ with random weights \label{step:distance_matrix} \Comment{Common Step 2}
\State \Return $D$
\end{algorithmic}
\end{algorithm}

\subsubsection{Spanning Tree Matrices}

Tree matrices are generated from random spanning trees. Similar to the bipartite matrix, the adjacency matrix is converted into a weighted symmetric matrix with randomly assigned weights. The steps are as follows:

\begin{algorithm}[H]
\caption{Random Tree Matrix Generation.}\label{alg:random_tree}
\begin{algorithmic}[1]
\Input $n$
\Output $D$: An $n \times n$ random tree weighted symmetric matrix
\State Generate a random tree $T$ with $n$ nodes
\State \textbf{Call} Algorithm \ref{alg:bipartite_mtrx}, Steps \ref{step:adjacency} and \ref{step:distance_matrix} respectively with input $T$ 
\State \Return $D$
\end{algorithmic}
\end{algorithm}

\subsubsection{Planar Graph Matrices}

To generate a planar matrix, we first create a complete graph $G$ with $n$ nodes. Then, with a brute force technique, we remove random edges until we obtain a planar structure. We check whether the graph we produce is connected, if not, we add a random edge. In the last stage, we perform the common steps, which are producing the adjacency matrix and weighted symmetric matrix respectively from $G$.

\begin{algorithm}[H]
\caption{Generate Planar Matrix.}\label{alg:planar_matrix}
\begin{algorithmic}[1]
\Input $n$
\Output $D$: An $n \times n$ planar weighted symmetric matrix
\State Generate a complete graph $G$ with $n$ nodes
\State Shuffle the edges of $G$
\State Set $is\_planar \gets \text{False}$
\While {not $is\_planar$}
    \State Remove an edge from $G$ 
    \State Check if $G$ is planar 
\EndWhile
\If {$G$ is not connected}
    \State Add an edge between two random nodes in $G$ 
\EndIf
\State \textbf{Call} Algorithm \ref{alg:bipartite_mtrx}, Steps \ref{step:adjacency} and \ref{step:distance_matrix} respectively with input $G$ 
\State \Return $D$
\end{algorithmic}
\end{algorithm}

\subsubsection{Chordal Graph Matrices}

 First, we generate a random tree $T$ consisting of $n$ nodes. In order to guarantee that the graph is chordal, we iterate over each node in T. For any node that has at least two neighboring nodes, we randomly choose two neighbors, denoted as $v$ and $w$, and add an edge between $v$ and $w$ to form a chord. This procedure preserves the chordal structure of the graph. The latter steps consist of transforming the adjacency matrix into a weighted symmetric matrix, a process that applies to all types of matrices previously mentioned.

\begin{algorithm}[H]
\caption{Generate Chordal Matrix.}\label{alg:chordal_matrix}
\begin{algorithmic}[1]
\Input $n$ 
\Output $D$: An $n \times n$ chordal weighted symmetric matrix
\State Generate a random tree $T$
\For {each node in $T$}
    \State Get the list of neighbors of the current node 
    \If {the number of neighbors is at least 2}
        \State Randomly select two neighbors $v$ and $w$ from $neighbors$
        \State Add an edge between $v$ and $w$ to create a chord
    \EndIf
\EndFor
\State \textbf{Call} Algorithm \ref{alg:bipartite_mtrx}, Steps \ref{step:adjacency} and \ref{step:distance_matrix} respectively with input $T$ 
\State \Return $D$
\end{algorithmic}
\end{algorithm}

\subsubsection{Odd-Cycle-Rich Matrices}

To generate a non-bipartite odd-cycle-rich matrix, we first construct a connected random graph and then explicitly enforce the presence of at least one triangle. This guarantees that the graph contains an odd cycle and is therefore non-bipartite. The resulting adjacency matrix is then converted into a weighted symmetric matrix. This family is particularly useful for generating structures in which SDP-based relaxations may dominate 
LP relaxations.

\subsubsection{Diagonal Matrices}

Next, we describe the generation of various types of $n \times n$ diagonal matrices used in our study. These matrices include those with ordered and randomly placed \(1\)s and \(-1\)s, as well as diagonal matrices with random values. In all cases, the resulting matrix \(D\) is a diagonal matrix formed as \texttt{\(D = \text{diag}(v)\)}, where \(v\) is a vector defined according to the specific method.

In this section, we produce diagonal matrices that vary in order and magnitude of their eigenvalues. The number of negative eigenvalues for these matrices, $n'$, can be specified or randomly determined. First, we create ordered diagonal matrices with eigenvalues $-1$ and $1$, where $n'$ values are $-1$ followed by \(n - n'\) values of $1$. In the second method, we give a random order to $-1$s and $1$s. And finally, we produce diagonal matrices with eigenvalues containing random numbers in random positions, provided that ($n'$) of them are negative. The aim here is to enhance the diversity of the matrices used in our tests.

\begin{algorithm}[H]
\caption{Diagonal Matrix Generation.}\label{alg:unified_diag_matrix}
\begin{algorithmic}[1]
\Input $n$, $type$ (one of \{`ordered\_ones', `random\_ones', `random\_randnums'\}), \(n'\) (number of negative eigenvalues, optional)
\Output $D$: An $n \times n$ diagonal matrix based on the specified type
\If {\(n'\) is not specified}
    \State Set \(n'\) to a random integer between 0 and $n$ (inclusive)
\EndIf
\State Initialize vector $v$ of length $n$
\If {$type$ is `ordered\_ones'}
    \State Set the first \(n'\) elements of $v$ to -1, and the remaining $(n - n')$ elements to 1
\ElsIf {$type$ is `random\_ones'}
    \State Set all elements of $v$ to 1
    \State Select \(n'\) random positions from the range $0$ to $n-1$ and store them in $neg\_list$
    \For {each index $idx$ in $neg\_list$}
        \State Set the element of the vector $v$ at position $idx$ to -1
    \EndFor
\ElsIf {$type$ is `random\_randnums'}
    \State Generate $n$ random values from the range 1 to 10 and store them in $v$
    \State Select \(n'\) random positions from the range $0$ to $n-1$ and store them in $neg\_list$
    \For {each index $idx$ in $neg\_list$}
        \State Set the element of the vector $v$ at position $idx$ to the negative of its current value
    \EndFor
\EndIf
\State Form a diagonal matrix $D$ from vector $v$
\State \Return $D$
\end{algorithmic}
\end{algorithm}


\subsubsection{Conversion of Adjacency to Weighted Symmetric Matrices}

For bipartite, tree, planar, and chordal matrices, the adjacency matrices are converted to weighted symmetric matrices. Each edge in the adjacency matrix is replaced with a random integer, ensuring symmetry and adding variability to the weights. The method is as follows:
\[
D[i,j] = \begin{cases} 
R_{ij}, & \text{if } A[i,j] = 1 \\
0, & \text{otherwise,}
\end{cases}
\]
Here, \( R_{ij} \) is a randomly selected integer $\sim \unif(-10, 10) $. 

            

\subsection{Random Matrices}\label{sec:random}

This part of the study describes the generation of various types of $n \times n$ random matrices used in our study. These include random symmetric matrices, matrices with specific eigenvalue properties, and random symmetric positive definite matrices. 

\subsubsection{Random Symmetric Matrix}

To generate random symmetric matrices, we began by creating a matrix \(A\) with elements drawn randomly $\sim \unif(-5, 5) $. 

\begin{algorithm}[H]
\caption{Random Symmetric Matrix Generation.}\label{alg:random_sym_mtrx}
\begin{algorithmic}[1]
\Input $n$
\Output $M$: An $n \times n$ random symmetric matrix
\State Generate a random symmetric $n \times n$ matrix $M$ with elements sampled from $ \unif(-5, 5) $
\State \Return $M$
\end{algorithmic}
\end{algorithm}

\subsubsection{Random Symmetric Positive Definite Matrix}

Random symmetric $n \times n$ positive definite matrices are generated using a standard procedure that ensures positive definiteness. The \texttt{scikit-learn} library \cite{scikit-learn}, a powerful tool for machine learning in Python, provides the function \texttt{make\_spd\_matrix} to generate symmetric positive definite (SPD) matrices.

\subsubsection{Random Matrices with Specified Eigenvalues}

Random matrices with specified eigenvalues are generated using a procedure involving QR decomposition. For matrices with eigenvalues of \(\pm 1\), first we generate random matrix, and perform its QR decomposition to obtain an orthogonal matrix $Q$. Then, we form a diagonal matrix \(\Lambda\) with entries of  \(1\) and  \(-1\), and the final matrix $M$ is constructed as  \(M = Q \Lambda Q^T\), representing its eigendecomposition. In the case of matrices with random eigenvalues, the process follows the same initial steps, but the diagonal matrix \(\Lambda\) is composed of random numbers selected from a specified range rather than fixed entries of \(1\) and \(-1\).

\begin{algorithm}[H]
\caption{Random Matrix with Specified Eigenvalues.}\label{alg:random_mtrx_specified_eigvals}
\begin{algorithmic}[1]
\Input $n$, $type$ (either `ones' or `random'), $n'$ (optional)
\Output $M$: Random $n \times n$ matrix with specified eigenvalues
\If {$n'$ is not specified} \label{step:check_specified}
    \State Randomly select $n'$ from 0 to $n$
\EndIf
\State Generate a random $n \times n$ matrix $rand\_mtrx$ with entries as continuous random numbers $\sim \unif(0, 9) $
\State  Generate an orthogonal matrix $Q$ \label{step:qr_decomp}
\State Initialize vector $v$ of length $n$
\If {$type$ is `ones'}
    \State Fill $v$ with ones
    \State Randomly select $n'$ positions in $v$ and set them to -1
\ElsIf {$type$ is `random'}
    \State Generate $v$ with random integers between 1 and 10
    \State Randomly select $n'$ positions in $v$ and set them to -1
\EndIf
\State Form a diagonal matrix $\Lambda$ using the vector $v$ \label{step:diag_matrix}
\State Construct the matrix $M$ from $\Lambda$ and $Q$  \label{step:construct_matrix}
\State \Return $M$
\end{algorithmic}
\end{algorithm}


\revise{
\section{Additional Computational Characteristics of Synthetic Instances}
\label{app:synthetic_computational_characteristics}

This appendix reports the detailed computational results for the
synthetic instances. The global optimization study uses 500 Synthetic
and 200 Synthetic$'$ instances, with a Gurobi time limit of 300 seconds
per instance and a relative optimality tolerance of $10^{-4}$.
Runtime and node statistics are calculated over instances for which
Gurobi certified optimality within this tolerance.

\begin{table}[htbp]
\centering
\caption{Global optimization outcomes and computational effort. Runtime and node statistics use certified-optimal instances only.}
\label{tab:synthetic_difficulty}
\begin{tabular}{lrr}
\toprule
Measure & Synthetic & Synthetic$'$ \\
\midrule
\# Distinct QCQPs & 500 & 200 \\
\# Proven optimal & 458 & 144 \\
Time-limit terminations & 42 & 56 \\
Median runtime, optimal cases (s) & 0.003 & 0.092 \\
Maximum runtime, optimal cases (s) & 291.476 & 294.377 \\
Median B\&B nodes, optimal cases & 1 & 21 \\
\bottomrule
\end{tabular}
\end{table}

For the time-limited instances in Table~\ref{tab:synthetic_timelimit_gaps},
we report the remaining Gurobi relative optimality gap. For a
minimization problem, it is calculated as
\begin{equation}
g_{\mathrm{opt}} =
\frac{ U-L }{\max\{|U|,10^{-10}\}},
\end{equation}
where $U$ is the best feasible incumbent objective value and $L$ is
the best global lower bound returned by Gurobi at termination.
The gap is reported only when both values are finite.

\begin{table}[htbp]
\centering
\caption{Remaining Gurobi relative optimality gaps for instances reaching the 300-second time limit (percent).}
\label{tab:synthetic_timelimit_gaps}
\begin{tabular}{lrr}
\toprule
Measure & Synthetic & Synthetic$'$ \\
\midrule
Time-limited instances & 42 & 56 \\
Median gap (\%) & 1.923 & 8.233 \\
90th percentile gap (\%) & 397.302 & 1872.082 \\
\bottomrule
\end{tabular}
\end{table}

\begin{table}[htbp]
\centering
\caption{Median normalized relaxation gaps to certified global optima. The number of instances with finite relaxation bounds is shown in parentheses.}
\label{tab:synthetic_gaps}
\begin{tabular}{lcc}
\toprule
Relaxation & Synthetic & Synthetic$'$ \\
\midrule
LP & 0.0406 (458) & 0.1033 (144) \\
SDP & $2.70\times10^{-8}$ (233) & 0.0246 (113) \\
SDP$'$ & $3.37\times10^{-7}$ (458) & 0.00185 (144) \\
\bottomrule
\end{tabular}
\end{table}

Table~\ref{tab:synthetic_root_distribution} reports root-node
relaxation comparisons over the complete datasets: 320,500 Synthetic
instances and 6,068 Synthetic$'$ instances.
Win and tie percentages are calculated over all comparisons, whereas
the median absolute relative difference uses finite relaxation bounds only.

\begin{table}[htbp]
\centering
\caption{Root-node relative bound differences over the complete synthetic datasets.}
\label{tab:synthetic_root_distribution}
\scriptsize
\setlength{\tabcolsep}{3pt}
\begin{tabular}{llrrrr}
\toprule
Dataset & Comparison & Median $|\delta_R|$ & Tie (\%) & SDP wins (\%) & LP wins (\%) \\
\midrule
Synthetic & LP--SDP & 0.3864 & 6.68 & 39.72 & 53.60 \\
Synthetic & LP--SDP$'$ & 0.0275 & 37.51 & 48.78 & 13.71 \\
Synthetic$'$ & LP--SDP & 0.2502 & 11.68 & 44.69 & 43.63 \\
Synthetic$'$ & LP--SDP$'$ & 0.1136 & 26.30 & 64.82 & 8.89 \\
\bottomrule
\end{tabular}
\end{table}}

\revise{
\section{Machine Learning Reproducibility Details}\label{app:hyperparameters}

The model configurations reported below were used consistently in the root node and spatial B\&B experiments. Hyperparameters were chosen manually rather than through a systematic grid search. The tables summarize the settings of the classical machine learning models and GNN, together with the data splitting and preprocessing procedures.

\begin{table}[h]
\centering
\small
\caption{Hyperparameters of the classical machine learning models used in the root node and spatial B\&B experiments. Classification and regression variants use the same settings unless noted otherwise.}
\label{tab:appendix_classical_hyperparameters}
\begin{tabular}{|p{0.2\linewidth}|p{0.7\linewidth}|}
\hline
\textbf{Model} & \textbf{Hyperparameter settings} \\
\hline
Extreme Gradient Boosting (classification and regression) & 400 boosting trees; learning rate 0.08; maximum tree depth 6; 80\% of training observations sampled per tree; 80\% of predictors sampled per tree; L1 penalty 0.1; L2 penalty 1.0; random seed 42. Classification uses logistic loss. \\
\hline
Gradient Boosting (classification and regression) & 100 boosting trees; learning rate 0.1; maximum depth of each regression tree 3; all training observations and predictors used per boosting iteration; random seed 42. Other settings follow the software defaults. \\
\hline
Random Forest (classification) & 400 decision trees; the number of candidate predictors considered at each split equals the square root of the total number of predictors; class weights inversely proportional to class frequencies; random seed 42. \\
\hline
Random Forest (regression) & 400 decision trees; all predictors considered at each split; bootstrap sampling enabled; random seed 42. \\
\hline
Data preprocessing & Missing numeric predictor values replaced with the corresponding training-set medians; the fitted imputation rules applied unchanged to the test data. No predictor standardization was applied to these models. \\
\hline
\end{tabular}
\end{table}

\begin{table}[h]
\centering
\small
\caption{Graph neural network architecture and training settings used in both the root node and spatial B\&B experiments.}
\label{tab:appendix_gnn_hyperparameters}
\begin{tabular}{|p{0.2\linewidth}|p{0.7\linewidth}|}
\hline
\textbf{Component} & \textbf{Specification} \\
\hline
Graph representation & Tripartite graph with variable nodes, quadratic-term nodes, and constraint nodes; five input features per variable node, five per quadratic-term node, and two per constraint node. Connections encode variable--quadratic-term, variable--constraint, and quadratic-term--constraint relationships in both directions. \\
\hline
Message passing & Three message-passing layers with 64 hidden units; permutation-invariant representation. Mean pooling produces a graph-level representation for each node type. \\
\hline
Prediction outputs & A classification output for relaxation selection, with an auxiliary prediction of the relative objective difference; a regression output for predicting the relative difference. \\
\hline
Optimization & Adam optimizer; 10 training epochs; learning rate 0.001; weight decay 0.00001; dropout probability 0.01; training batch size 8 and evaluation batch size 16. \\
\hline
Loss functions & Classification: weighted binary cross-entropy with an auxiliary squared-error loss (coefficient 0.25). Regression: weighted squared-error loss with an auxiliary decision-classification loss (coefficient 0.50; decision-logit scale 20). Sample-weight coefficient 0.30. \\
\hline
Feature preprocessing & Node and edge features standardized using means and standard deviations computed on training graphs only; features with a standard deviation below $10^{-8}$ use a divisor of 1. No handcrafted structural features or early-tree probe features were provided as GNN inputs. \\
\hline
\end{tabular}
\end{table}

\begin{table}[h]
\centering
\small
\caption{Dataset partitioning and experimental comparison protocol. These procedures determine which instances were used for fitting and evaluating the models.}
\label{tab:appendix_splits}
\begin{tabular}{|p{0.2\linewidth}|p{0.7\linewidth}|}
\hline
\textbf{Setting} & \textbf{Procedure} \\
\hline
Training/test split & Random split of the combined MINLPLib and QPLIB dataset (random seed: 42). \\
\hline
Training on synthetic instances & Models were trained on the Synthetic$'$ dataset and evaluated on the combined MINLPLib--QPLIB benchmark. No instances from the combined evaluation benchmark were included in synthetic training. \\
\hline
Consistency across learning tasks & The same eligible benchmark and instance partitions were used for classification and regression.
\\

\hline
\end{tabular}
\end{table}
}

\section{GNN Prediction Diagnostics}
\label{app:gnn_diagnostics}

\revise{To investigate the root node GNN results, Table~\ref{tab:gnn_diagnostics} reports additional classification and regression diagnostics on benchmark instances.}

\begin{table}[htbp]
\centering
\caption{Root-node GNN prediction diagnostics on benchmark instances.}
\label{tab:gnn_diagnostics}
\scriptsize
\resizebox{\textwidth}{!}{%
\begin{tabular}{llrrrrrrr}
\hline
& & & & \multicolumn{2}{c}{GNN-C} & \multicolumn{2}{c}{GNN-R} & \\
\cline{5-6}\cline{7-8}
Comparison & Training & $n$ & True LP (\%) &
Pred. LP (\%) & Acc. (\%) &
Pred. LP (\%) & r-Acc. (\%) & Majority (\%) \\
\hline
LP--SDP & Combined & 101 & 73.27 & 81.19 & 92.08 & 80.20 & 89.11 & 73.27 \\
LP--SDP & Synthetic$'$ & 278 & 80.58 & 35.25 & 48.20 & 83.09 & 66.55 & 80.58 \\
LP--SDP$'$ & Combined & 101 & 56.44 & 67.33 & 87.13 & 64.36 & 85.15 & 56.44 \\
LP--SDP$'$ & Synthetic$'$ & 278 & 64.75 & 31.65 & 44.60 & 10.79 & 58.63 & 64.75 \\
\hline
\end{tabular}%
}
\end{table}

\revise{The synthetic-trained GNN-C models predict LP substantially less frequently than its true occurrence in the benchmark instances, resulting in accuracies below the majority-class baseline. GNN-R similarly exhibits lower r-accuracy under synthetic training. These findings are consistent with weaker cross-dataset generalization of both GNN models.}

\bibliographystyle{spmpsci}      

\bibliography{Library}           

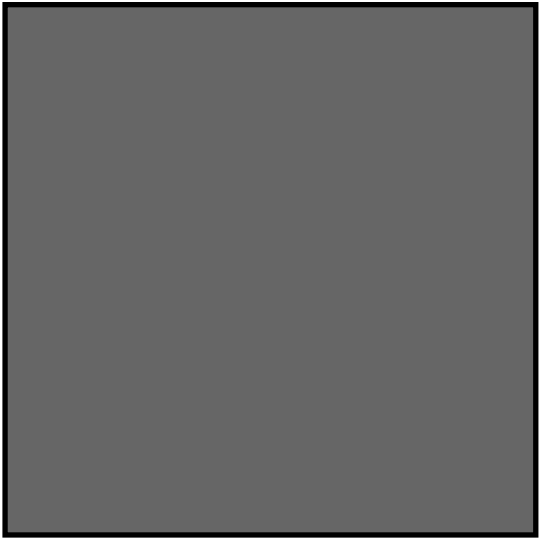
\end{document}